\documentclass[number,citesort,dvips]{arxbj}
\usepackage{upgreek}
\usepackage{graphicx}

\aid{0}
\volume{18}
\issue{1}
\pubyear{2012}
\firstpage{137}
\lastpage{176}
\doi{10.3150/10-BEJ332}

\makeatletter
\newtheorem{theorem}{Theorem}
\newremark{remark}{Remark}
\newtheorem{lemma}{Lemma}
\makeatother

\begin{document}
\begin{frontmatter}

\title{On asymptotically optimal wavelet estimation of trend functions under
long-range dependence}

\runtitle{On wavelet trend estimation under long-range dependence}

\begin{aug}
\author{\fnms{Jan} \snm{Beran}\corref{}\thanksref{e1}\ead[label=e1,mark]{Jan.Beran@uni-konstanz.de}}
\and
\author{\fnms{Yevgen} \snm{Shumeyko}\thanksref{e2}\ead[label=e2,mark]{Yevgen.Shumeyko@uni.konstanz.de}}

\runauthor{J. Beran and Y. Shumeyko}
\address{Department of Mathematics and Statistics,
University of Konstanz, Universit\"{a}tsstra\ss e 10, \mbox{D-78464}
Konstanz, Germany.\\ \printead{e1,e2}}
\end{aug}

% HISTORY:
\received{\smonth{7} \syear{2009}}
\revised{\smonth{9} \syear{2010}}

% ABSTRACT
\begin{abstract}
We consider data-adaptive wavelet estimation of a trend function in
a time series model with strongly dependent Gaussian residuals.
Asymptotic expressions for the optimal mean integrated squared error
and corresponding optimal smoothing and resolution parameters are
derived. Due to adaptation to the properties of the underlying trend
function, the approach shows very good performance for smooth trend
functions while remaining competitive with minimax wavelet
estimation for functions with discontinuities. Simulations
illustrate the asymptotic results and finite-sample behavior.
\end{abstract}

% KEYWORDS

\end{frontmatter}

%s1 ###
\section{Introduction}\vspace*{-3pt}

Suppose that we observe time series data of the form
\begin{equation}\label{1}
Y_{i}=g(t_{i})+\xi_{i},\qquad i=1,2,\dots,n,
\end{equation}
with $t_{i}=i/n$, $g\in L^{2}( [0,1]) $ and $\xi_{i}$ a Gaussian
zero-mean second order stationary process with long-range dependence.
Here, long-range dependence is characterized by
\begin{equation}\label{5}
\gamma(k)=E( \xi_{i}\xi_{i+k}) \mathop{\sim}_{k\rightarrow\infty}
C_{\gamma}|k|^{-\alpha}
\end{equation}
for some constants $\alpha\in(0,1)$ and $C_{\gamma}>0,$ where `$\sim$'
means that the ratio of the two sides converges to $1$. For the
spectral density
$f(\lambda)=(2\uppi)^{-1}\sum\gamma(k)\exp(-\mathrm{i}k\lambda),$ this
corresponds to a pole at the origin of the form $C_{f}|\lambda| ^{\alpha-1}$ for a suitable constant $C_{f}$.

Nonparametric estimation of $g$ in this context has been studied
extensively in the last two decades, including kernel smoothing
(Hall and Hart \cite{r32}, Cs\"{o}rg\"{o} and Mielniczuk \cite{r18,r19}, Ray and Tsay \cite{r42}, Robinson \cite{r43}, Beran and Feng \cite{r11,r12}),
local polynomial\vadjust{\goodbreak} estimation (Beran and Feng \cite{r13}, Beran \textit{et al.} \cite{r14}) and wavelet
thresholding (Wang \cite{r48}, Johnstone and Silverman \cite{r37}). For nonparametric quantile
estimation in long-memory processes, see also Ghosh \textit{et al.} \cite{r28} and
Ghosh and Draghicescu \cite{r29,r30}. In this paper, we take a~closer look at optimal
wavelet estimation of $g$. Wang~\cite{r48} and Johnstone and Silverman~\cite{r37} derived optimal
minimax rates within general function spaces and Gaussian
long-memory residuals. In particular, the minimax threshold
$\sigma\sqrt{2\log n}$ turns out to achieve the minimax rate even
under long memory. However, for some practical applications, the minimax
approach may  be too pessimistic. It may, for instance, be
known a~priori that $g$ or some derivatives of $g$ are piecewise
continuous. Li and Xiao \cite{r39} therefore considered data-adaptive selection
of resolution levels. They derived an asymptotic expansion for the
mean integrated squared error (MISE) under the assumptions that $g$
is piecewise smooth and the resolution levels used for the
estimation are chosen according to certain asymptotic rules
(formulated in terms of the parameters $J$ and $q$, as defined
below). The rate of the MISE achieved this way turns out to be the
same as for minimax rules. No further justification for the specific
choice of $J$ and $q$ is given, however, and no optimality result is
derived. We refer to Remark \ref{remark2} below for further
discussion on Li and Xiao \cite{r39}.\looseness=-1

In this paper, the aim is to obtain concrete data-adaptive rules for
optimal estimation of~$g$. In a first step, it is shown that for
functions with continuous derivatives, the rate given in Li and Xiao \cite{r39}
can be achieved without thresholding by choosing optimal values of~$J$
and $q.$ In a second step, exact constants for the MISE and
asymptotic formulas for the optimal choice of $J$ and $q$ are
derived. These results are comparable to results on optimal bandwidth
selection in kernel smoothing (Gasser and M\"{u}ller~\cite{r27}, Hall and Hart~\cite{r32},
Beran and Feng \cite{r11,r13}). In a third step, additional higher resolution
levels combined with thresholding are added in order to include the
possibility of discontinuities. The resulting estimator shows very
good performance for smooth trend functions (comparable to optimal
kernel estimators) while remaining competitive with (and even
superior to) minimax wavelet estimation for functions with
discontinuous derivatives.

For literature on trend estimation by wavelet thresholding in the
case of i.i.d.~or weakly dependent residuals, see, for example,
Donoho and Johnstone \cite{r22,r25,r23}, Donoho \textit{et al.}~\cite{r24}, Daubechies \cite{r21}, Brillinger \cite{r15,r16},
Abramovich \textit{et al.} \cite{r5}, Nason \cite{r40}, Johnstone and Silverman \cite{r37}, Johnstone \cite{r36}, Percival and Walden \cite{r41},
Vidakovic \cite{r47}, Hall and Patil~\cite{r33,r34,r35}, Sachs and Macgibbon \cite{r45} and Truong and Patil \cite{r44}. Apart
from Johnstone and Silverman \cite{r37} and Wang \cite{r48}, wavelet trend estimation in the
long-memory case has also been considered by Yang \cite{r51} for random
design models.

The paper is organized as follows. Basic definitions are introduced
in Section \ref{section2}. The main results are given in Section
\ref{section3}. A simulation study in Section \ref{section4}
illustrates the results. Concluding remarks are given in Section
\ref{section5}. Proofs can be found in the \hyperref[app]{Appendix}.

%s2 ###
\vspace*{-3pt}\section{Basic definitions}\vspace*{-3pt}\label{section2}

Let $\phi(t)$ and $\psi(t)$ be the father and mother wavelets,
respectively, with compact support $[0,N]$ for some $N\in\mathbb{N}$
and such that
\begin{eqnarray}
\label{intphi2}\int_{0}^{N}\phi(t)\,\mathrm{d}t&=&\int_{0}^{N}\phi^{2}(t)\,\mathrm{d}t=\int_{0}^{N}\psi^{2}(t)\,\mathrm{d}t=1,
\\
\label{11}\psi(0)&=&\psi(N)=0\vadjust{\goodbreak}
\end{eqnarray}
and, for any $J\geq0$, the system $\{\phi_{Jk},\psi_{jk},k\in\mathbb{Z}%
,j\geq0\}$ with
\[
\psi_{jk}(t)=N^{1/2}2^{(J+j)/2}\psi(N2^{J+j}t-k),\qquad
\phi_{Jk}(t)=N^{1/2}2^{J/2}\phi(N2^{J}t-k)
\]
is an orthonormal basis in $L^{2}(\mathbb{R})$. Note that for the
sake of generality, the support of~$\phi$ and $\psi$ is chosen to be
$[0,N]$ instead of $[0,1]$. This way, it is possible to choose from
a~larger variety of wavelet generating functions satisfying
(\ref{intphi2}) (see Daubechies \cite{r21}, Cohen \textit{et al.}~\cite{r17}). Throughout the paper,
$m_{\psi}\in\mathbb{N}$ will denote the number of vanishing moments
of~$\psi$, that is,
\begin{equation}\label{9}
\int_{0}^{N}t^{k}\psi(t)\,\mathrm{d}t=0,\qquad k=0,1,\ldots,m_{\psi}-1,
\end{equation}
and
\begin{equation}\label{10}
\int_{0}^{N}t^{m_{\psi}}\psi(t)\,\mathrm{d}t=\nu_{m_{\psi}}\neq0.
\end{equation}
For every function $g\in L^{2}([0,1])$ and every $J\geq0$, we have
the orthogonal wavelet expansion
\begin{equation}\label{3}
g(t)=\sum_{k=-N+1}^{N2^{J}-1}s_{Jk}\phi_{Jk}(t)+\sum_{j=0}^{\infty}\sum
_{k=-N+1}^{N2^{J+j}-1}d_{jk}\psi_{jk}(t),
\end{equation}
where
\[
s_{Jk}=\int_{0}^{1}g(t)\phi_{Jk}(t)\,\mathrm{d}t,\qquad
d_{jk}=\int_{0}^{1}g(t)\psi _{jk}(t)\,\mathrm{d}t
\]
are the wavelet coefficients of the function $g$. A (hard)
thresholding wavelet estimator of~$g$ is defined by
\begin{equation}\label{4}
\hat{g}(t)=\sum_{k=-N+1}^{N2^{J}-1}\hat{s}_{Jk}\phi_{Jk}(t)+\sum_{j=0}^{q}%
\sum_{k=-N+1}^{N2^{J+j}-1}\hat{d}_{jk}I(|\hat{d}_{jk}|>\delta_{j})%
\psi_{jk}(t),
\end{equation}
where $J$, $q$ and $\delta_{j}$ denote the decomposition level,
smoothing
parameter and threshold, respectively, and the wavelet coefficients $\hat {s}%
_{Jk}$ and $\hat{d}_{jk}$ are given by
\[
\hat{s}_{Jk}=\frac{1}{n}\sum_{i=1}^{n}Y_{i}\phi_{Jk}(t_{i})\quad\mbox{and}\quad\hat{d}_{jk}=%
\frac{1}{n}\sum_{i=1}^{n}Y_{i}\psi_{jk}(t_{i});
\]
see, for example, Donoho and Johnstone \cite{r22,r23}, Abramovich
\textit{et al.} \cite{r5}. For estimates without thresholding (i.e.,
$\delta_{j}\equiv0$), see also  Johnstone and Silverman \cite{r37} and
Nason \cite{r40}, Brillinger \cite{r15,r16}, among others.

%s3 ###
\section{Main results}\label{section3}
\vspace*{-3pt}

In the context of long-memory errors, an explicit asymptotic
expansion for the MISE is given in Li and Xiao \cite{r39} under specific
assumptions on the decomposition level $J$ and the smoothing
parameter $q$. The question of how to choose $J$ and $q$
optimally is not investigated. The following theorem establishes
the optimal convergence rate of the MISE when minimizing with
respect to $J$, $q$ and $\{\delta_j\}$.

 In what follows, $\phi$
and $\psi$ will  be assumed either to be piecewise differentiable or
to satisfy a uniform H\"{o}lder condition with exponent $1/2$, that is,
\begin{equation}\label{Hoelder}
|\psi(x)-\psi(y)| \le C|x-y|^{1/2} \qquad \forall x,y \in [0,N] .
\end{equation}
Daubechies (\cite{r21}, Chapter 6)  provides examples of wavelets
satisfying these conditions. Moreover, throughout this paper,
$2^{J}=\mathrm{o}(n)$ to ensure that $\hat g$ includes resolution
levels lower than the distance between successive time points. This
assumption is needed for the consistency of $\hat g$, as discussed
below.

\begin{theorem}\label{theorema1}
Suppose that $g\in C^{r}[0,1]$, the support
$\operatorname{supp}(g^{(r)})=\{t\in\lbrack0,1]\dvtx g^{(r)}(t)\neq0\}$
has positive Lebesgue measure, the process $\xi_{i}$ is Gaussian with
covariance structure \textup{(\ref{5})} and $\psi$ is such that
$m_{\psi}=r$. Then, minimizing the MISE with respect to $J$, $q$ and
$\{\delta_j\}$ yields the optimal order
\begin{equation}
\mathit{MISE}_{\mathrm{opt}}=\mathrm{O}\bigl(n^{-2r\alpha/(2r+\alpha)}\bigr).
\end{equation}
\end{theorem}

Theorem \ref{theorema1} is of limited practical use since only rate
optimality is established.  Theorem \ref{Theorem2}
will show that the rate obtained in Li and Xiao \cite{r39} can be achieved without
thresholding {by minimizing the MISE }with respect to $J$ and $q$.
In order to apply the result to observed data, optimal constants
need to be derived. This question is addressed in Theorems
\ref{Theorem2} and \ref{theorema3} below. The following constants
will be needed:
\begin{eqnarray}
\label{C_phi}
C_{\phi}^{2}
&=&
C_{\gamma}\int_{0}^{N}\int_{0}^{N}|x-y|^{-\alpha}\phi(x)\phi(y)\,\mathrm{d}x\,\mathrm{d}y,
\\[-2pt]
\label{C_psi}
C_{\psi}^{2}
&=&
C_{\gamma}\int_{0}^{N}\int_{0}^{N}|x-y|^{-\alpha}\psi(x)\psi(y)\,\mathrm{d}x\,\mathrm{d}y,
\\[-2pt]
C^{\ast}\bigl(r,\alpha,\psi,g^{(r)}\bigr)
&=&
\frac{1}{2r+\alpha}\log_{2}\biggl[ \frac {\int_{0}^{N}\nu_{r}^{2}(g^{(r)}(t))^2\,\mathrm{d}t}{C_{\psi}^{2}(r!)^{2}}\biggr]-\log_{2}N,\nonumber
\\[-2pt]
\label{DeltaPsi}\Delta_n(g,C_{\psi})&=&\frac{\alpha}{2r+\alpha}\log_2n+C^{\ast}\bigl(r,\alpha,\psi,g^{(r)}\bigr)\nonumber
\\[-9pt]\\[-9pt]
&&{}-\biggl\lfloor\frac{\alpha}{2r+\alpha}\log_2n+C^{\ast}\bigl(r,\alpha,\psi,g^{(r)}\bigr)\biggr\rfloor,\nonumber
\end{eqnarray}
where $\lfloor x \rfloor$ denotes the largest integer less than or
equal to $x$,
\begin{eqnarray}
A_{1}(r,\alpha,\psi) & =&\biggl( \frac{2^{2r\Delta_n(g,C_{\psi})}}{2^{2r}-1}+\frac{2^{\alpha(1-\Delta_n(g,C_{\psi}))}}{%
2^{\alpha}-1}\biggr) ( C_{\psi}^{2}) ^{2r/(2r+\alpha)},\nonumber
\\
A_{2}\bigl(r,\alpha,\psi,g^{(r)}\bigr) & =&
\biggl( \frac{\nu_{r}^{2}}{(r!)^{2}}\int
_{0}^{1}\bigl(g^{(r)}(t)\bigr)^2\,\mathrm{d}t\biggr) ^{\alpha/(2r+\alpha)},\nonumber
\\
\nu_{r}&=&\int t^{r}\psi(t)\,\mathrm{d}t,\nonumber
\\
C^{\ast}\bigl(r,\alpha,\phi,g^{(r)}\bigr)&=&\frac{1}{2r+\alpha}\log_{2}\biggl[ \frac{%
\int_{0}^{1}\nu_{r}^{2}(g^{(r)}(t))^2\,\mathrm{d}t}{C_{\phi}^{2}(2^{\alpha
}-1)(r!)^{2}}\biggr] -\log_{2}N,\nonumber
\\
\label{DeltaPhi}\Delta_n(g,C_{\phi})&=&\frac{\alpha}{2r+\alpha}\log_2n+C^{\ast}\bigl(r,\alpha,\phi,g^{(r)}\bigr)\nonumber
\\[-8pt]\\[-8pt]
&&{}-\biggl\lfloor\frac{\alpha}{2r+\alpha}\log_2n+C^{\ast}\bigl(r,\alpha,\phi,g^{(r)}\bigr)\biggr\rfloor,\nonumber
\\
A_{3}(r,\alpha,\phi)&=&\biggl(
\frac{2^{2r\Delta_n(g,C_{\phi})}}{2^{2r}-1}+\frac{2^{\alpha(1-\Delta_n(g,C_{\psi}))}}{2^{\alpha
}-1}\biggr) \bigl( C_{\phi}^{2}(2^{\alpha}-1)\bigr) ^{2r/(2r+\alpha)}.\nonumber
\end{eqnarray}

For the case where no thresholding is used, exact asymptotic
expressions for the MISE and an optimal solution can be given as
follows.

\begin{theorem}
\label{Theorem2}
Under the assumptions of Theorem \ref{theorema1} and thresholds%
\[
\delta_{j}=0\qquad(0\leq j\leq q),
\]
the following holds.
\begin{enumerate}[(ii)]
\item[(i)] If $(2^{\alpha}-1)C_{\phi}^{2}>C_{\psi}^{2}$, then the asymptotic
MISE is minimized by the smoothing parameter
\begin{equation}\label{Bestq1}
q^{\ast}=\biggl\lfloor\frac{\alpha}{2r+\alpha}\log_{2}n+
C^{\ast}\bigl(r,\alpha,\psi,g^{(r)}\bigr)\biggr\rfloor-J^{\ast}
\end{equation}
with decomposition levels $J^{\ast}$ satisfying
$2^{J^{\ast}}=\mathrm{o}(n^{\alpha/(2r+\alpha)})$. The optimal MISE is of
the form
\begin{equation}\label{16}
\mathit{MISE}=A_{1}(r,\alpha,\psi)A_{2}\bigl(r,\alpha,\psi,g^{(r)}\bigr)\cdot n^{-2r\alpha/(2r+\alpha)}+\mathrm{o}\bigl( n^{-2r\alpha/(2r+\alpha)}\bigr) .
\end{equation}

Moreover, if $\Delta_n(g,C_{\psi})=0$, then
\[
q^{\ast}=\biggl\lfloor\frac{\alpha}{2r+\alpha}\log_{2}n+
C^{\ast}\bigl(r,\alpha,\psi,g^{(r)}\bigr)\biggr\rfloor-J^{\ast}-1
\]
(with $J^{\ast}$ as before) also minimizes the $\mathit{MISE}$.
\item[(ii)] If $(2^{\alpha}-1)C_{\phi}^{2}<C_{\psi}^{2}$, then minimizing
the asymptotic MISE with respect to $J$ and $q$ yields
\begin{equation}\label{bestj2}
J^{\ast}=\biggl\lfloor\frac{\alpha}{2r+\alpha}\log_{2}n+C^{\ast}\bigl(r,\alpha,\phi
,g^{(r)}\bigr)\biggr\rfloor+1
\end{equation}
and
\begin{equation}
\hat{g}(t)=\sum_{k=-N+1}^{N2^{J}-1}\hat{s}_{Jk}\phi_{Jk}(t)
\end{equation}
with $J=J^{\ast}$. The optimal MISE is of the form
\begin{equation}\label{19}
\mathit{MISE}=A_{3}(r,\alpha,\phi)A_{2}\bigl(r,\alpha,\psi,g^{(r)}\bigr)\cdot n^{-2r\alpha/(2r+\alpha)}+\mathrm{o}\bigl( n^{-2r\alpha/(2r+\alpha)}\bigr).
\end{equation}
Moreover, if $\Delta_n(g,C_{\phi})=0$, then
\[
J^{\ast}=\biggl\lfloor\frac{\alpha}{2r+\alpha}\log_{2}n+C^{\ast}\bigl(r,\alpha,\phi,g^{(r)}\bigr)\biggr\rfloor
\]
also minimizes the $\mathit{MISE}$.
\end{enumerate}
\end{theorem}

If higher resolution levels beyond those used in Theorem
\ref{Theorem2} are included together with thresholding, then the
values of the MISE given in (\ref{16}) and (\ref{19})
can be attained even if~$g^{(r)}$ does not exist everywhere and is
only piecewise continuous.

\begin{theorem}\label{theorema3}
Suppose {that} {$g^{(r)}$}{exists on $[0,1]$ except for {at most a}
finite number of points and, where it exists, it is piecewise
continuous and bounded}. Furthermore, assume that~%
$\operatorname{supp}(g^{(r)})$ has positive Lebesgue measure,
$m_{\psi}=r$ and the process $\xi_{i}$ is Gaussian and such that
(\ref{5}) holds. The following then hold:
\begin{enumerate}[(ii)]
\item[(i)] if $(2^{\alpha}-1)C_{\phi}^{2}>C_{\psi}^{2}$, $J$ is such
that $2^{J}=\mathrm{o}(n^{\alpha/(2r+\alpha)})$, $q=\lfloor
\log_{2}n\rfloor -J $, $q^{\ast}$ is defined by (\ref{Bestq1}) and
$\delta_{j}$ is such that for $0\leq j\leq q^{\ast}$,
\begin{equation}
\delta_j=0,
\end{equation}
and for $q^{\ast}<j\leq q,$
\begin{equation}
2^{J+j}\delta_{j}^{2}\rightarrow 0,\mbox{ }2^{(J+j)(2r+1)}\delta_{j}^{2}\rightarrow\infty,\qquad \delta_j^2\ge
\frac{4 \mathrm{e} C_{\psi}^2N^{-1+\alpha}(\ln n)^2
}{n^{\alpha}2^{(J+j)(1-\alpha)}},
\end{equation}
then equation (\ref{16}) holds;
\item[(ii)] if $(2^{\alpha}-1)C_{\phi}^{2}<C_{\psi}^{2}$, $J=J^{\ast}$ with $%
J^{\ast}$ defined by (\ref{bestj2}), $q=\lfloor \log_{2}n\rfloor -J$ and
$\delta_{j}$ is such that
\begin{eqnarray}
2^{J+j}\delta_{j}^{2}&\rightarrow&0,\qquad2^{(J+j)(2r+1)}\delta_{j}^{2}%
\rightarrow\infty,\nonumber
\\[-8pt]\\[-8pt]
 \delta_j^2&\ge& \frac{4 \mathrm{e} C_{\psi}^2N^{-1+\alpha}(\ln
n)^2}{n^{\alpha}2^{(J+j)(1-\alpha)}}\qquad(0\leq j\leq q),\nonumber
\end{eqnarray}
then equation (\ref{19}) holds.
\end{enumerate}
\end{theorem}

\begin{remark}
Li and Xiao \cite{r39} derived an asymptotic expansion for the MISE
under the assumptions that $J,q\rightarrow\infty,$
$2^{J+j}\delta_{j}^{2}\rightarrow0$,
$2^{(2r+1)(J+j)}\delta_{j}^{2}\rightarrow\infty$ and $\delta_{j}^{2}$
are\vspace*{1pt} above a certain bound that\vadjust{\goodbreak} depends on $j,$ $n$, $g,$ $\alpha$ and
$J.$ The question of how to choose $J$, $q$ and
$\delta_{j}$ optimally is not considered. Here, a partial solution to
the optimality problem is given. Theorem \ref{Theorem2} provides
optimal values of $q$ and $J$, and a corresponding
formula for the optimal MISE, for estimators with no thresholding (i.e., $%
\delta_{j}\equiv0$). This result is obtained for $r$-times
continuously
differentiable trend functions. Thus, jumps and other irregularities in $g$
are excluded. In a second step, we therefore ask the question whether
the asymptotic formula for the optimal MISE can be extended to more
general functions. Theorem \ref{theorema3} shows that this is indeed
the case, in the sense that (essentially) $g$ does not need to be
differentiable everywhere. This includes, for instance, the
possibility of isolated jumps. Note that for a given $n$, $q=\lfloor
\log_{2}n\rfloor -J$ is the highest available resolution. By adding all
available higher resolution levels combined with thresholding, the same
formula for the MISE applies {as in Theorem \ref{Theorem2}}. The
intuitive reason for this is that isolated discontinuities are
`infinitesimally local' and can therefore be characterized best when
the finest possible levels of resolution are included. At very high
resolution, however, non-zero thresholds are needed in order to
distinguish deterministic jumps from noise. For functions where Theorem
\ref{Theorem2} applies, the optimal MISE in Theorem \ref{Theorem2} and
the MISE obtained in Theorem~\ref{theorema3} are the same.
\end{remark}

\begin{remark}\label{remark2}
The only quantity in (\ref{Bestq1}) and (\ref{bestj2}) that depends on
$n$ is $\alpha(2r+\alpha
)^{-1}\log_{2}n$. The constants $C^{\ast}(r,\alpha,\psi,g^{(r)})$ and $%
C^{\ast}(r,\alpha,\phi,g^{(r)})$ provide data-adaptive
adjustments to optimize the multiplicative constant in the MISE.
They can be decomposed into several terms with different meanings.
For instance,
\[
C^{\ast}\bigl(r,\alpha,\phi,g^{(r)}\bigr)=\frac{C_{1}^{\ast}+C_{2}^{\ast}+C_{3}^{%
\ast}}{2r+\alpha}+C_{4}^{\ast}
\]
with
\[
C_{1}^{\ast}=\log_{2}\int_{0}^{1}\bigl(g^{(r)}(t)\bigr)^2\,\mathrm{d}t
\]
reflecting the properties of $g$,
\[
C_{2}^{\ast}=\log_{2}\biggl( \frac{\nu_{r}}{r!}\biggr) ^{2}
\]
depending on the basis function $\psi$,
\[
C_{3}^{\ast}=-\log_{2}[ C_{\phi}^{2}(2^{\alpha}-1)]
\]
characterized by the basis function $\phi$ and the asymptotic
covariance structure (\ref{5}) of $\xi_{i}$, and
\[
C_{4}^{\ast}=-\log_{2}N
\]
defined by the length of the support of $\psi$ and $\phi$. Note that for $%
N=1$, $C_{4}^{\ast}=0$.
\end{remark}

\begin{remark}
The question of how far the MISE can be optimized further with respect
to freely adjustable thresholds is more difficult and is the subject of
current research. The same comment applies to the possibility of soft
thresholding.
It is worth mentioning here, however, that for some classes of functions, $%
\delta_{j}=0$ is indeed the best threshold. For instance, it can be
shown that if $g\in L^{2}[0,1]$ and $C<|g^{(r)}(\cdot)|\leq
C2^{r+\alpha/2}$ (almost everywhere) for some finite constant $C$, then
$\delta_{j}=0$ is asymptotically optimal. This includes, for example,
functions that can be represented (or approximated in an appropriate
sense) by piecewise $r$th order polynomials.\vspace*{2pt}
\end{remark}

\begin{remark}
The results in Li and Xiao \cite{r39} are derived for residuals of the form $\xi
_{i}=G(Z_{i}),$ where~$Z_{i}$ is a stationary Gaussian long-memory
process and the transformation~$G$ has Hermite rank $m_{G}$. For
simplicity of presentation, the results given here are only derived for
Gaussian processes. An extension to $\xi_{i}=G(Z_{i})$ would be
possible along the same lines.\vspace*{2pt}
\end{remark}

\begin{remark}
Asymptotic expressions for the MISE and formulas for optimal
bandwidth selection in kernel regression with long memory are given
in Hall and Hart \cite{r32}, Cs\"{o}rg\"{o} and Mielniczuk \cite{r18} and Beran and Feng \cite{r11,r13}, among others. Note,
however, that there,~$g^{(r)}$ has to be assumed to be continuous
instead of only piecewise continuous, and $r\geq2$. In that sense,
the applicability of kernel estimators (and also of local
polynomials) is more limited. This is illustrated in the simulation
study in the next section.\looseness=1\vspace*{2pt}
\end{remark}

\begin{remark}
In analogy to kernel estimation, the optimal rate of convergence of
wavelet estimates becomes faster the more derivatives of $g$ that exist.
However, the optimal MISE can only be achieved if the number of
vanishing moments of the mother wavelet~$\psi$ is equal to $r$. In
other words, the choice of an appropriate wavelet basis is
essential. This is analogous to kernel estimation where a kernel of
the appropriate order should be used (see, e.g., Gasser and M\"{u}ller \cite{r27}).
Consider, for instance, the case where only the first derivative of
$g$ exists (and is piecewise continuous), that is, $r=1$.
Then, for the wavelets estimator, the optimal order of the MISE is $\mathrm{O}(n^{-2\alpha/(2+\alpha)})$. In this case, we may use Haar wavelets
(for which $m_{\psi}=1$). In contrast to the wavelet estimator, the
usual asymptotic expansion for the MISE of kernel estimators does
not hold in this case. On the other hand, if $g$ is twice
continuously differentiable, then
the optimal rate achieved by kernel estimators is at least $\mathrm{O}(n^{-4\alpha/(4+\alpha)})$. If Haar wavelets are used, then, in spite of $r$
being equal to $2$,  the optimal rate of the wavelet estimator
cannot be better than $\mathrm{O}(n^{-2\alpha/(2+\alpha)})$ and is
thus slower than the rate achieved by kernel estimators. In order to
match the rate of kernel estimators, a wavelet basis with
$m_{\psi}=2$ vanishing moments has to be used.\vspace*{2pt}
\end{remark}

\begin{remark}
The optimal rate of convergence of the MISE is the same as the
minimax rate obtained by Wang \cite{r48} and Johnstone and Silverman \cite{r37}. However, for a
given function, the multiplicative constant in the asymptotic
expression of the MISE is essential. This is achieved here by data-adaptive choices of $q$ and $J$. The simulations in the next section
illustrate that the data-adaptive method tends to outperform the
minimax solution, provided that the assumptions of Theorems
\ref{Theorem2} or \ref{theorema3} hold.
\end{remark}

\begin{remark}
The best smoothing parameter and decomposition level depend on the
unknown parameters $\alpha$, $C_{\gamma}$ and the unknown $r$th
derivative of $g$. Based on Theorems \ref{Theorem2} and
\ref{theorema3}, an iterative data-adaptive algorithm along the lines
of Beran and Feng \cite{r12} can be designed. Essentially, the iteration consists
of a step where $g$ is estimated (using the best estimates of
relevant parameters available at that stage) and a step where $\alpha$, $%
C_{\gamma}$ and other quantities in the asymptotic MISE formula are
estimated. For the estimation of $C_{\gamma}$ and $\alpha$, see, for
instance, Yajima \cite{r49}, Fox and Taqqu \cite{r26}, Dahlhaus \cite{r20}, Giraitis and Surgailis \cite{r31}, Beran \cite{r8,r9}, Beran \textit{et al.} \cite{r10}, Abry and Veitch \cite{r6}. A detailed iterative algorithm is
currently being developed and will be presented elsewhere. {An
obvious choice for estimating $\alpha$ is to use an appropriate
wavelet-based method such as that described in Bardet \textit{et al.} \cite{r7}.} Note that
while the idea of the iteration is simple, a concrete implementation
is far from trivial (see Beran and Feng \cite{r12}). In particular, in the presence
of long-range dependence, small changes in the smoothing parameters
can lead to considerable changes in the estimate of the long-memory
parameter $\alpha,$ and vice versa.\vspace*{-6pt}
\end{remark}

%f1 ###
\begin{figure}[b]

\includegraphics{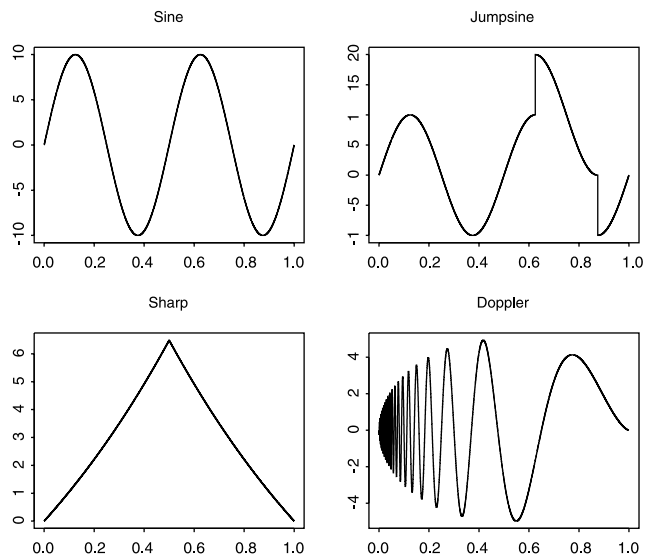}

\caption{Trend functions used
in the simulations: sine, JumpSine, ``sharp'' and Doppler.}\label{figure1}
\end{figure}

%s4 ###
\section{Simulations}
\vspace*{-3pt}\label{section4}

To study the potential benefits of data-adaptive wavelet estimation
as outlined above, a simulation study was carried out with four
different test
functions $g$ (Figure \ref{figure1}) and a~Gaussian $\operatorname{FARIMA}(0,d,0$) residual process $%
\xi_{i}$. Note that $\alpha=1-2d$. The test functions are:
\begin{itemize}
\item sine function: $g_{1}(t)=10\sin(4\uppi t);$\vadjust{\goodbreak}
\item JumpSine function: $g_{2}(t)=10\sin(4\uppi t)+\Delta\cdot I\{\frac{5}{8}%
<t<\frac{7}{8}\}$ $(\Delta > 0);$
\item ``sharp'' function: $g_{3}(t)=10[ \exp( tI\{ t<0.5\}
+(1-t)I\{t>0.5\}) -1]; $
\item Doppler function: $g_{4}(t)=10[ t(1-t)] ^{1/2}\sin[
2\uppi(1+0.05)/(t+0.05)]. $
\end{itemize}

The following methods are compared:
\begin{itemize}
\item Wavelet estimator with hard thresholding, $q$, $J$ as in Theorem \ref{theorema3} and
\[
\delta_j^2=\frac{4 \mathrm{e} C_{\psi}^2N^{-1+\alpha}(\ln n)^2
}{n^{\alpha}2^{(J+j)(1-\alpha)}}\qquad  (q^{\ast}<j\leq q).
\]
Note that for the first three functions, Theorem \ref{theorema3}(ii)
applies, whereas for the Doppler function, derivatives are not
bounded. Nevertheless, we carried out the simulations using a
modified version of $C^{\ast}$ (see the remarks at the end of this
section).
\item Wavelet estimator with soft thresholding defined by
\[
\operatorname{sign}(\hat d_{jk})(|\hat d_{jk}|-\lambda_{n})I\{|\hat
d_{jk}|>\lambda_{n}\}
\]
and minimax thresholds
\[
\lambda_{n}=(2\log n)^{1/2}
\]
(Johnstone and Silverman \cite{r37}).
\item Kernel estimator with rectangular kernel $K(x)=\frac{1}{2}I\{x\in
\lbrack-1,1]\}$ and asymptotically optimal bandwidth
\[
b_{\mathrm{opt}}=C_{\mathrm{opt}}n^{(2d-1)/(5-2d)},
\]
where
\begin{eqnarray*}
C_{\mathrm{opt}}&=&\biggl(\frac{9(1-2d)\beta(d)C_{f}}{I(g^{\prime\prime})}\biggr)
^{1/(5-2d)},
\\
\beta(d)&=&\frac{2^{2d}\Gamma(1-2d)\sin(\uppi d)}{d(2d+1)}
\end{eqnarray*}
(see, e.g., Hall and Hart \cite{r32}, Beran and Feng \cite{r11}).
\end{itemize}
%f2 ###
\begin{figure}

\includegraphics{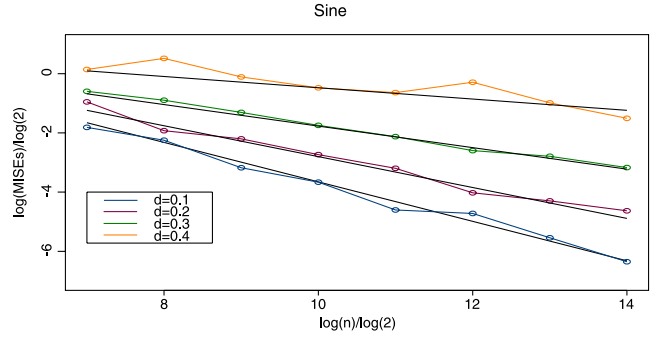}

\caption{Simulated values of
the mean
 integrated squared error, $\mathit{MISE}_{\mathrm{sim}}$, for different
values of the fractional
parameter $d$, plotted against the sample size ($n=2^{7},2^{8},\ldots,2^{13}$%
) on log--log scale (base 2 logarithms). The results are based on
$400$ simulations of model (\protect\ref{1}) with the sine trend function
and $\operatorname{FARIMA}(0,d,0$) residuals with $d=0.1, 0.2, 0.3,0.4.$ The
estimates are based on Theorem \protect\ref{theorema3} and wavelet basis
$s4$.}\label{figure2}
\vspace*{12pt}
\end{figure}
%t1 ###
\begin{table}
\caption{Logarithms
(base 2) of simulated values of the mean integrated squared
error, $\log_{2} \mathit{MISE}_{\mathrm{sim}},$ as a function of $n$ and the wavelet
bases s4, s6, s8 and s10, respectively. For comparison,
$\log_{2} \mathit{MISE}_{\mathrm{theor}}$ obtained from the asymptotic formulas in
Theorem \protect\ref{theorema3} is also given. The results are based on $400
$ simulations of a $\operatorname{FARIMA}(0,0.2,0)$ model with trend function $%
g_{1}(t)=10\sin(4\protect\uppi t)$}\label{table1}
\begin{tabular*}{\textwidth}{@{\extracolsep{\fill}}lllll@{}}
\hline
$n$ & Simulation `s4' & Theor. `s4' & Simulation `s6' & Theor. `s6' \\
\hline
\phantom{00$\,$}128 & 0.516420047 & 0.408553554 & 0.251744659 & 0.332459614\\
\phantom{00$\,$}256 & 0.263441364 & 0.294451230 & 0.214928924 &0.222321976 \\
\phantom{00$\,$}512 & 0.217604044 & 0.219171771 & 0.112951872& 0.149658234 \\
\phantom{0$\,$}1024 & 0.150284851 & 0.150545678 &0.110547951 & 0.101718042 \\
\phantom{0$\,$}2048 & 0.109213215 & 0.100879757& 0.079795806 & 0.070089311 \\
\phantom{0$\,$}4096 & 0.061483507 &0.068112469 & 0.049441935 & 0.049222131 \\
\phantom{0$\,$}8192 & 0.050871673& 0.046494121 & 0.030814609 & 0.035454926 \\
16$\,$384 &0.040330363 & 0.032231330 & 0.020141994 & 0.026371959 \\[6pt]
$n$ & Simulation `s8' & Theor.~`s8' & Simulation `s10' & Theor.~`s10'\\
\hline
\phantom{00$\,$}128 & 0.251744659 & 0.290131091 & 0.348379471 &0.251989178 \\
\phantom{00$\,$}256 & 0.214928924 & 0.193352318 & 0.20541786 &0.174618829 \\
\phantom{00$\,$}512 & 0.112951872 & 0.129502140 & 0.158692616& 0.123573436 \\
\phantom{0$\,$}1024 & 0.110547951 & 0.087376732 &0.074319167 & 0.089896035 \\
\phantom{0$\,$}2048 & 0.079795806 & 0.059584328& 0.061712354 & 0.065326166 \\
\phantom{0$\,$}4096 & 0.049441935 &0.041248179 & 0.030175723 & 0.043107368 \\
\phantom{0$\,$}8192 & 0.030814609& 0.029150833 & 0.027662929 & 0.028448428 \\
16$\,$384 &0.020141994 & 0.021169561 & 0.020361623 & 0.018777135 \\
\hline
\end{tabular*}%
\end{table}
\textit{Sine}: Figure \ref{figure2} shows reasonably good agreement
between the simulated and theoretical MISE of the adaptive wavelet
estimator with basis s4. Here, s4, s6$,\ldots$ denote Daubechies'
wavelets with $2,3,\dots$ vanishing moments, respectively (see
Daubechies \cite{r21}).
Table \ref{table1} illustrates the effect of using different basis functions for the case $%
d=0.2$. Irrespective of the wavelet basis (s4, s6, s8 or s10), the agreement between the simulated MISE and the
theoretical formula is already very good  for $n=256.$ However, since
$g$ is infinitely continuously differentiable, the MISE can be
reduced by using very smooth basis functions. This explains why the
performance of s4 is considerably worse compared with s6, s8 and s10.
Table \ref{table2} shows that, as expected, the mean squared error
increases with increasing long memory (see also Figure
\ref{figure2}). A comparison between minimax wavelet thresholding,
the data-adaptive wavelet estimator and kernel smoothing is given in
Figures \ref{figure3} and \ref{figure4}. Since the sine function is
well behaved, optimal kernel estimation is expected to perform well.
The kernel estimator does indeed outperform the minimax procedure.
In contrast, the MISE of the data-adaptive wavelet method is
comparable to optimal kernel estimation. A typical sample path and
the corresponding estimated trend functions are plotted in Figure
\ref{figure5}. The minimax rule leads to a rather erratic function
near local minima and maxima, whereas this is not the case for the
other two methods.

%t2 ###
\begin{table}
\caption{Simulated values of
the MISE for different sample sizes and values of $d$. The results
are based on $400$ simulations of model (\protect\ref{1}) with
$\operatorname{FARIMA}(0,d,0)$ residuals, the sine trend function~$g_{1}$ and the
wavelet estimator based on Theorem \protect\ref{theorema3} with wavelet
basis s4}\label{table2}
\begin{tabular*}{\textwidth}{@{\extracolsep{\fill}}lllll@{}}
\hline
$n$ & $d=0.1$ & $d=0.2$ & $d=0.3$ & $d=0.4$ \\
\hline
\phantom{$\,$00}128 & 0.284521469 & 0.516420047 &0.661787865 & 1.104194018 \\
\phantom{$\,$00}256 & 0.210694474 & 0.263441364& 0.537558642 & 1.42979724 \\
\phantom{$\,$00}512 & 0.110584545 & 0.217604044& 0.403889173 & 0.927229839 \\
\phantom{$\,$0}1024 & 0.078905169 &0.150284851 & 0.29832426 & 0.717419015 \\
\phantom{$\,$0}2048 & 0.041133887& 0.109213215 & 0.228981208 & 0.64283222 \\
\phantom{$\,$0}4096 &0.037871696 & 0.061483507 & 0.165045782 & 0.818104781 \\
\phantom{$\,$0}8192& 0.021438157 & 0.050871673 & 0.1444763 & 0.505236717 \\
16$\,$384 & 0.012234701 & 0.040330363 & 0.11107171 & 0.351823994 \\
\hline
\end{tabular*}%
\end{table}

%f3 ###
\begin{figure}[b]

\includegraphics{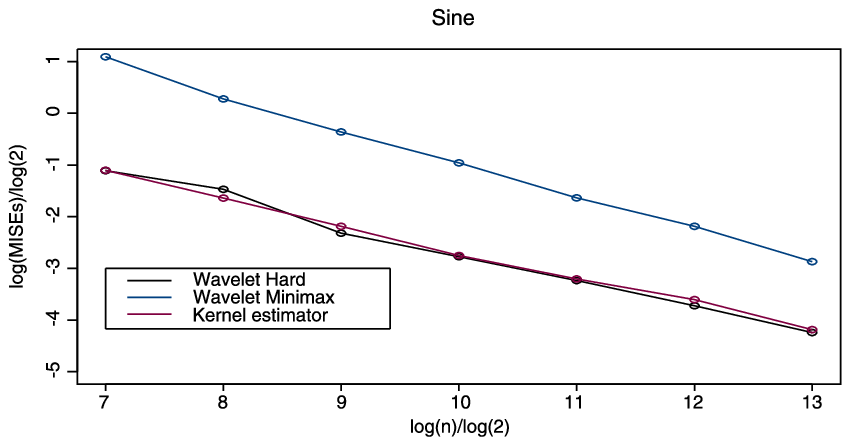}

\caption{Simulated values of $\log_{2}
\mathit{MISE}_{\mathrm{sim}}$ plotted against $\log n$
($n=2^{7},2^{8},\ldots,2^{13}$) for trend estimates obtained by kernel
smoothing, minimax soft threshold wavelet estimation and data-adaptive
hard threshold wavelet estimation obtained from Theorem \protect\ref{theorema3}
(both with basis~s4). The results are based on $400$ simulations of
model (\protect\ref{1}) with the sine trend function and $\operatorname{FARIMA}(0,0.2,0)$
residuals.}\label{figure3}
\end{figure}

\textit{Jumpsine}: The simulated and asymptotic MISE for the
Jumpsine function are compared in Table \ref{table3} for $d=0.2$ and jump sizes $%
\Delta=0.1,0.5,1,10,20$ and $50$. The agreement between the
asymptotic and simulated MISE is reasonably good, in particular for
small and very large values of $\Delta.$ Figure \ref{figure6}a shows
a typical sample path with $d=0.3$ and fits obtained by the three
methods. Figure \ref{figure6}b shows that, as expected from Theorem~\ref{theorema3}(ii), almost all non-zero coefficients belong to the
father wavelet. The mother wavelet functions are useful for modeling the
two jumps. Due to thresholding, almost all coefficients are
eliminated except those near $t=5/8$ and $7/8$. Similar
results were obtained for other values of $d. $ In comparison, the
data-adaptive wavelet method shows the best performance (Figures
\ref{figure7} and \ref{figure8}), although the difference between the
two wavelet methods is smaller under strong long memory. As
expected, kernel estimation cannot compete with the wavelet
approach.

%f4 ###
\begin{figure}

\includegraphics{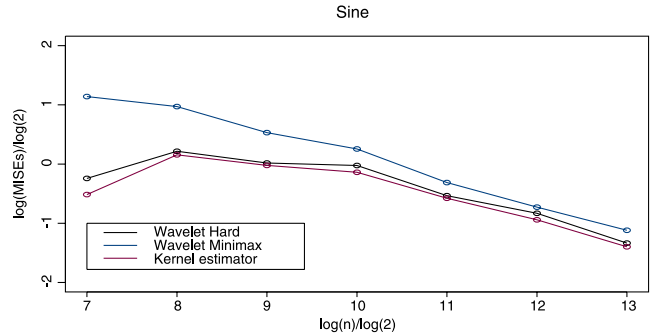}

\caption{Simulated values of $\log_{2}
\mathit{MISE}_{\mathrm{sim}}$ plotted against $\log n$
($n=2^{7},2^{8},\ldots,2^{13}$) for trend estimates obtained by kernel
smoothing, minimax soft threshold wavelet estimation and data-adaptive
hard threshold wavelet estimation obtained from Theorem \protect\ref{theorema3}
(both with basis~s4). The results are based on $400$ simulations of
model (\protect\ref{1}) with the sine trend function and
$\operatorname{FARIMA}(0,0.4,0)$ residuals.}\label{figure4}
\vspace*{2pt}
\end{figure}

%f5 ###
\begin{figure}[b]
\vspace*{2pt}
\includegraphics{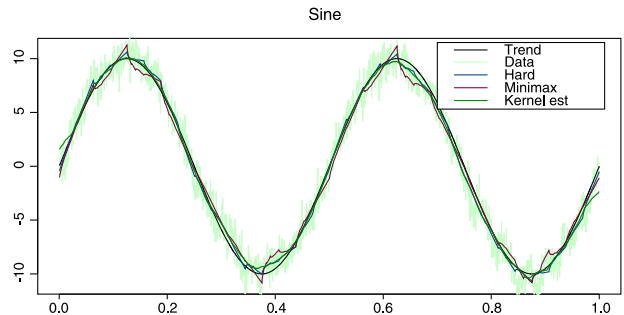}

\caption{Simulated data with sine function plus $%
\operatorname{FARIMA}(0,0.3,0)$ process, and trend estimates obtained by optimal
kernel smoothing, minimax soft thresholding wavelet estimation and
data-adaptive hard threshold wavelet estimation according to Theorem
\protect\ref{theorema3} (both with basis s4).}\label{figure5}
\end{figure}

%t3 ###
\begin{table}
\tablewidth=300pt
\caption{$\mathit{MISE}_{\mathrm{sim}}/\mathit{MISE}_{\mathrm{theor}}$ for the JumpSine
function and FARIMA(0, 0.2, 0) residuals, in dependence on the jump size $%
\Delta.$ The results are based on $400$ simulations and a
thresholding estimate according to Theorem \protect\ref{theorema3}, with
wavelet basis $s4$}\label{table3}
\begin{tabular*}{300pt}{@{\extracolsep{\fill}}llll@{}}
\hline
$\Delta$ & $n=2048$ & $n=4096$ & $n=8192$ \\
\hline
 \phantom{0}0.1 &1.02984365 & 1.000066053 & 0.996328962 \\
 \phantom{0}0.5 & 1.044736472 &1.007194657 & 1.004583086 \\
 \phantom{0}1 & 1.10352021 & 1.120497921 &1.096100157 \\
 10 & 1.635074083 & 1.690840646 & 1.563330038\\
 20 & 1.301618649 & 1.234763386 & 1.207770083 \\
 50& 1.222581848 & 1.21888936 & 1.115174282 \\
\hline
\end{tabular*}%
\end{table}

%f6 ###
\begin{figure}[b]

\includegraphics{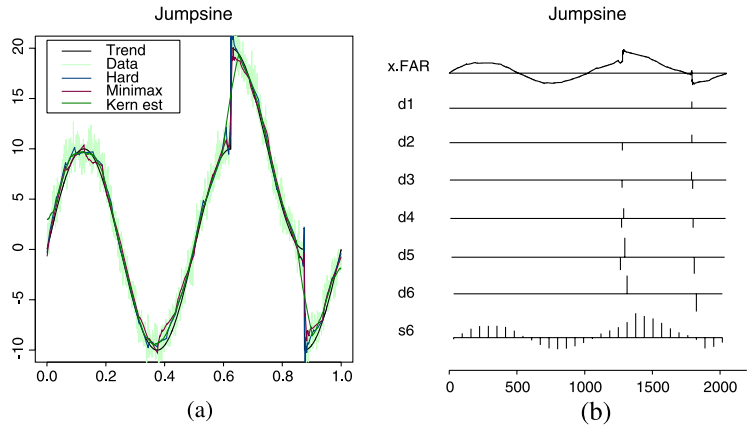}

\caption{Simulated data
(a) with JumpSine function plus $\operatorname{FARIMA}(0,0.3,0)$ process,
and trend estimates obtained by kernel smoothing, minimax soft
threshold wavelet estimation and data-adaptive hard threshold
wavelet estimation obtained from Theorem \protect\ref{theorema3} (both with
basis s4); (b) shows the coefficients of the data-adaptive
wavelet estimate.}\label{figure6}
\end{figure}

\textit{Sharp}: In distinct contrast to the JumpSine function, for the
sharp function, the performance of the kernel estimator is
comparable to the data-adaptive wavelet method (Figures
\ref{figure9} and \ref{figure10}), at least when the criterion is
the MISE. With respect to the visual fit, as exemplified by Figure
\ref{figure11}, the kernel method leads to oversmoothing of the edge
in the middle.

%f7 ###
\begin{figure}

\includegraphics{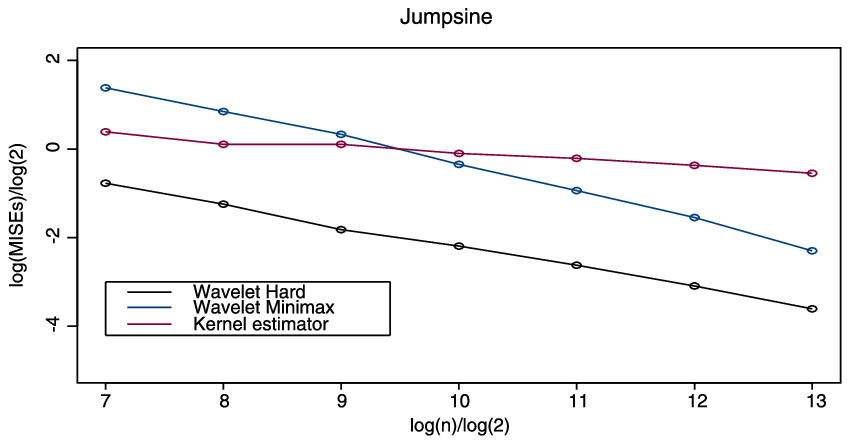}

\caption{Simulated values of
$\log_{2} \mathit{MISE}_{\mathrm{sim}}$ plotted against $\log n$
($n=2^{7},2^{8},\ldots,2^{13}$) for trend estimates obtained by
kernel smoothing, minimax soft threshold wavelet estimation and data-adaptive hard threshold wavelet estimation obtained from Theorem
\protect\ref{theorema3} (both with basis s4). The results are based on $400$
simulations of model (\protect\ref{1}) with the JumpSine trend function and
$\operatorname{FARIMA}(0,0.2,0)$ residuals.}\label{figure7}
\end{figure}

%f8 ###
\begin{figure}[b]

\includegraphics{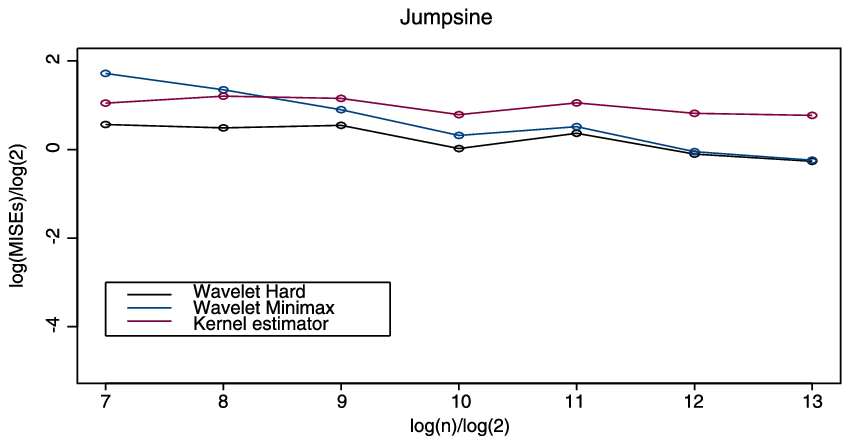}

\caption{Simulated values of
$\log_{2} \mathit{MISE}_{\mathrm{sim}}$ plotted against $\log n$
($n=2^{7},2^{8},\ldots,2^{13}$) for trend estimates obtained by
kernel smoothing, minimax soft threshold wavelet estimation and data-adaptive hard threshold wavelet estimation obtained from Theorem
\protect\ref{theorema3} (both with basis s4). The results are based on $400$
simulations of model (\protect\ref{1}) with the JumpSine trend function and
$\operatorname{FARIMA}(0,0.4,0)$ residuals.}\label{figure8}
\end{figure}

\textit{Doppler}: For the Doppler function, Theorem \ref{theorema3} is not applicable and $%
J^{\ast}$ in equation~(\ref{bestj2}) is not well defined.
Nevertheless, it is interesting to see how well hard thresholding may
work with a slight
modification of (\ref{bestj2}). Specifically, consider%
\[
\tilde{J}^{\ast}=\biggl\lfloor\frac{\alpha}{2r+\alpha}\log_{2}n+\tilde{C}%
^{\ast}\bigl(r,\alpha,\psi,\phi,g^{(r)}\bigr)\biggr\rfloor+1,
\]
where
\[
\tilde{C}^{\ast}\bigl(r,\alpha,\psi,\phi,g^{(r)}\bigr)=\frac{1}{2r+\alpha}\log _{2}%
\biggl[ \frac{\int_{0.1}^{0.95}\nu_{r}^{2}(g^{(r)}(t))^2\,\mathrm{d}t}{%
C_{\phi}^{2}(2^{\alpha}-1)(r!)^{2}}\biggr] -\log_{2}N.
\]
Note that the only change compared to $C^{\ast}$ consists of
bounding the
integration limits away from $0$ and $1$. For moderate long memory with $%
d=0.2,$ the data-adaptive wavelet estimator still turns out to be
the best (Figure \ref{figure12}). For strong long memory with $d=0.4$,
the minimax approach appears to be slightly better for very long
series (Figure \ref{figure13}). The relatively good performance of
the minimax approach is expected because, in contrast to the data-adaptive estimator, the coarser levels of resolution are not favored
a priori. This way, it is easier to catch the increasingly fast
oscillations toward the left of the timescale. As expected, the
kernel method does not work well. A typical example is shown in
Figure~\ref{figure14}.

%f9 ###
\begin{figure}

\includegraphics{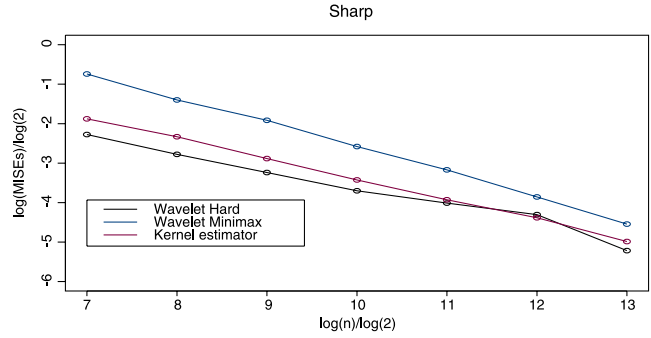}

\caption{Simulated values of
$\log_{2} \mathit{MISE}_{\mathrm{sim}}$ plotted against $\log n$
($n=2^{7},2^{8},\ldots,2^{13}$) for trend estimates obtained by
kernel smoothing, minimax soft threshold wavelet estimation and data-adaptive hard threshold wavelet estimation obtained from Theorem
\protect\ref{theorema3} (both with basis s4). The results are based on $400$
simulations of model (\protect\ref{1}) with the ``sharp'' trend function and
$\operatorname{FARIMA}(0,0.2,0)$ residuals.}\label{figure9}
\vspace*{3pt}
\end{figure}

%f10 ###
\begin{figure}[b]
\vspace*{3pt}
\includegraphics{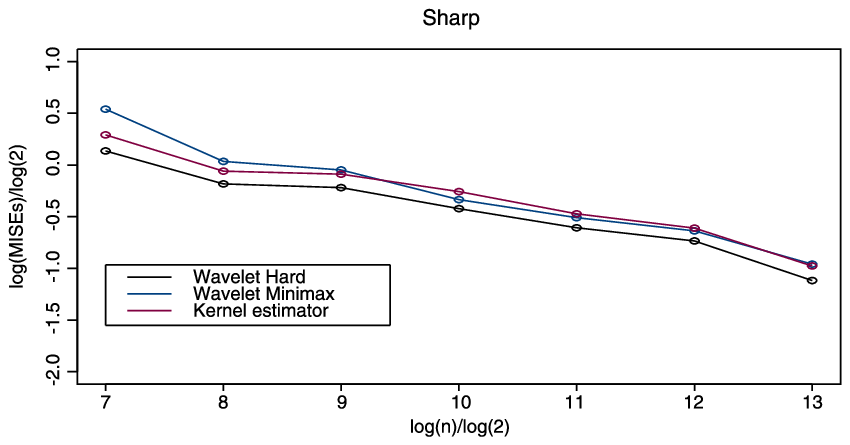}

\caption{Simulated values of
$\log_{2} \mathit{MISE}_{\mathrm{sim}}$ plotted against $\log n$
($n=2^{7},2^{8},\ldots,2^{13}$) for trend estimates obtained by
kernel smoothing, minimax soft threshold wavelet estimation and data-adaptive hard threshold wavelet estimation obtained from Theorem
\protect\ref{theorema3} (both with basis s4). The results are based on $400$
simulations of model (\protect\ref{1}) with the ``sharp'' trend function and
$\operatorname{FARIMA}(0,0.4,0)$ residuals.}\label{figure10}
\end{figure}

%f11 ###
\begin{figure}

\includegraphics{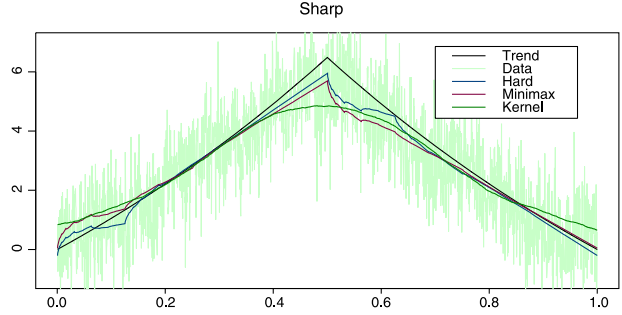}

\caption{Simulated data with ``sharp'' function plus $%
\operatorname{FARIMA}(0,0.3,0)$ process, and trend estimates obtained by kernel
smoothing, minimax soft threshold wavelet estimation and data-adaptive hard threshold wavelet estimation obtained from Theorem
\protect\ref{theorema3} (both with basis s4).}\label{figure11}
\end{figure}

%f12 ###
\begin{figure}[b]
\includegraphics{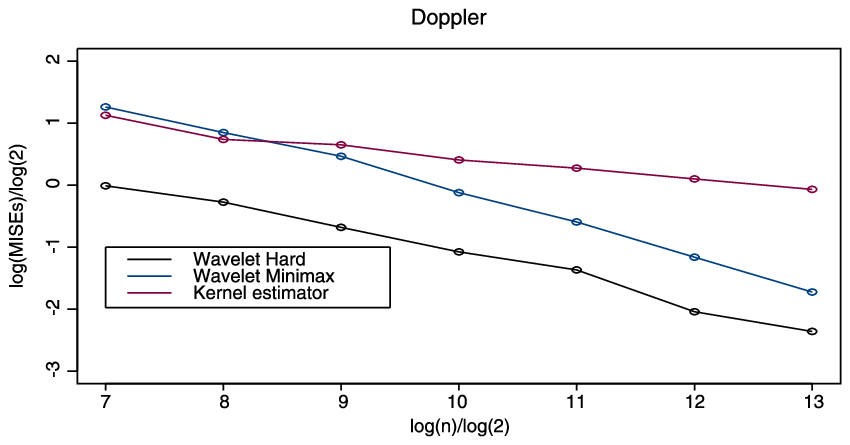}

\caption{Simulated values of
$\log_{2} \mathit{MISE}_{\mathrm{sim}}$ plotted against $\log n$
($n=2^{7},2^{8},\ldots,2^{13}$) for trend estimates obtained by
kernel smoothing, minimax soft threshold wavelet estimation and data-adaptive hard threshold wavelet estimation with $J=\tilde{J}^{*}$
and thresholds $\delta_{i}$ as in Theorem~\protect\ref{theorema3}(ii) (both with basis s4). The results are based on $400$
simulations of model (\protect\ref{1}) with the Doppler trend function and
$\operatorname{FARIMA}(0,0.2,0)$ residuals.}\label{figure12}
\end{figure}

%f13 ###
\begin{figure}

\includegraphics{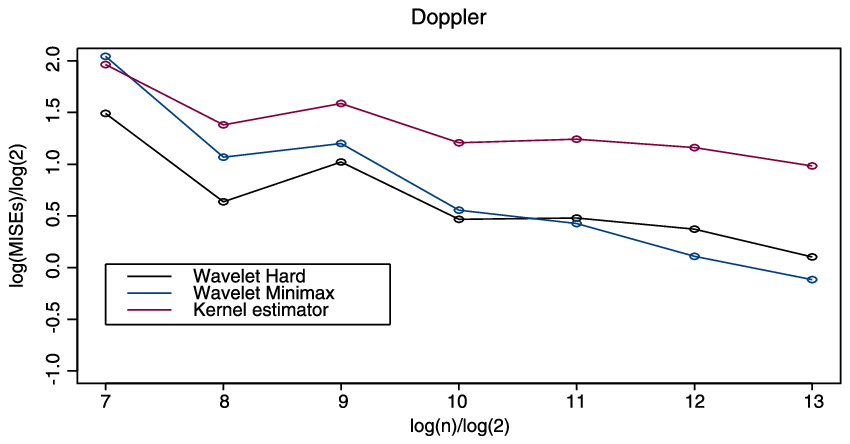}

\caption{Simulated values of
$\log_{2} \mathit{MISE}_{\mathrm{sim}}$ plotted against $\log n$
($n=2^{7},2^{8},\ldots,2^{13}$) for trend estimates obtained by
kernel smoothing, minimax soft threshold wavelet estimation and data-adaptive hard threshold wavelet estimation with $J=\tilde{J}^{*}$
and thresholds $\delta_{i}$ as in Theorem~\protect\ref{theorema3}(ii) (both with basis s4). The results are based on $400$
simulations of model (\protect\ref{1}) with the Doppler trend function and
$\operatorname{FARIMA}(0,0.4,0)$ residuals.}\label{figure13}
\end{figure}

%f14 ###
\begin{figure}[b]

\includegraphics{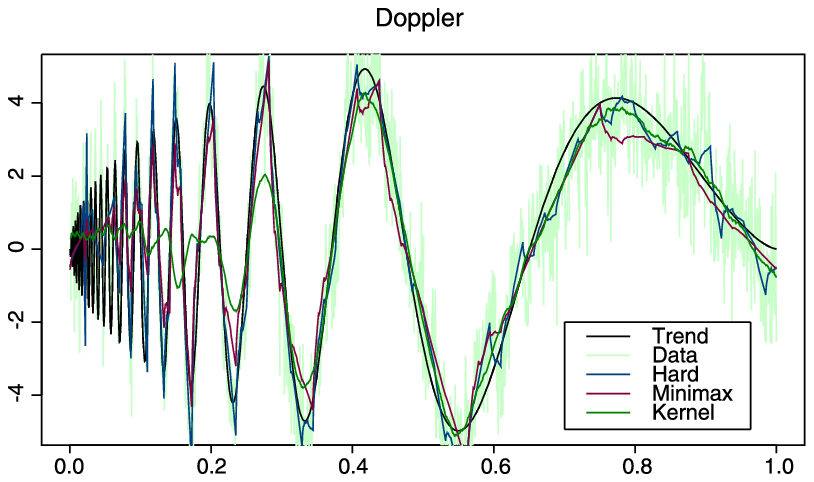}

\caption{Simulated data with
the Doppler function plus $\operatorname{FARIMA}(0,0.3,0)$ process, and trend
estimates obtained by kernel smoothing, minimax soft threshold
wavelet estimation and data-adaptive hard
threshold wavelet estimation with $J=\tilde{J}^{*}$ and thresholds $\delta_{i}$ as in Theorem~\protect\ref{theorema3}(ii) (both with basis s4).}\label{figure14}
\end{figure}

%s5 ###
\section{Concluding remarks}
\vspace*{2pt}
\label{section5}

In this paper, an approach to data-adaptive wavelet estimation of
trend functions for long-memory time series models is proposed. The
estimator can be understood as a~combination of two components: a
smoothing component consisting of a certain number of lower
resolution levels where no thresholding is applied and a higher
resolution component filtered by thresholding. The first component
leads to good performance for smooth functions, whereas the second
component is useful for modeling discontinuities. An open problem
worth pursuing in future research is the question of how much more may
be gained by further optimization with respect to fully flexible
thresholds $\delta_{j}$.

\begin{appendix}\label{app}
\renewcommand{\theequation}{\arabic{equation}}
\setcounter{equation}{22}
\section*{Appendix: Proofs}

In the proofs of Theorems \ref{theorema1}, \ref{Theorem2} and
\ref{theorema3}, $\phi$ and $\psi$ will be assumed to be piecewise
differentiable. Analogous results (apart from some expressions in
the remainder terms) can be obtained even if $\phi^{\prime}$ and
$\psi^{\prime}$ do not exist anywhere, provided that both functions
$\phi$ and $\psi$ satisfy a uniform H\"{o}lder condition with
exponent $1/2$ (see (\ref{Hoelder})). The proofs are analogous, with
the difference that instead of the rectangle rule (\ref{8}), the mean
value theorem is applied.

\begin{pf*}{Proof of Theorem \ref{theorema1}}
Let
\begin{equation}\label{6}
\mathit{MISE}=E\biggl[
\int_{0}^{1}\bigl(g(t)-\hat{g}(t)\bigr)^{2}\,\mathrm{d}t\biggr]\vadjust{\goodbreak}
\end{equation}
denote the mean integrated square error. Combining (\ref{6}) with
(\ref{3}) and (\ref{4}), we have
\begin{eqnarray*}
&&\mathit{MISE}=E\Biggl\{ \int_{0}^{1}\Biggl[ \sum_{k=-N+1}^{N2^{J}-1}(s_{Jk}-\hat{s}%
_{Jk})\phi_{Jk}(t)
\\
 &&\hphantom{\mathit{MISE}=E\Biggl\{ \int_{0}^{1}\Biggl[}{}+\sum_{j=0}^{q}\sum_{k=-N+1}^{N2^{J+j}-1}\bigl( d_{jk}-\hat{d}_{jk}I(|%
\hat{d}_{jk}|>\delta_{j})\bigr) \psi_{jk}(t)
\\
&&\hphantom{\mathit{MISE}=E\Biggl\{ \int_{0}^{1}\Biggl[}{}+ \sum
_{j=q+1}^{\infty}\sum_{k=-N+1}^{N2^{J+j}-1}d_{jk}\psi_{jk}\Biggr]
^{2}\,\mathrm{d}t\Biggr\} .
\end{eqnarray*}
Orthonormality of the basis in $L^{2}(\mathbb{R})$ implies
that
\begin{eqnarray}\label{MISE}
\mathit{MISE}
&=&
E\Biggl\{ \sum_{k=-N+1}^{N2^{J}-1}[\hat{s}_{Jk}-s_{Jk}]^{2}\Biggr\}\nonumber
\\
&&{}+
E\Biggl\{ \sum_{j=0}^{q}\sum_{k=-N+1}^{N2^{J+j}-1}[ \hat {d}_{jk}I(|\hat{d}_{jk}|>\delta_{j})-d_{jk}] ^{2}\Biggr\}+\sum_{j=q+1}^{\infty}\sum_{k=-N+1}^{N2^{J+j}-1}d_{jk}^{2}\nonumber
\\
&=&
\sum_{k=-N+1}^{N2^{J}-1}[ E( \hat{s}_{Jk}) -s_{Jk}] ^{2}+\sum _{k=-N+1}^{N2^{J}-1}E\{ [\hat{s}_{Jk}-E(\hat{s}_{Jk})] ^{2}\}
\\
&&{}+
\sum_{j=0}^{q}\sum_{k=-N+1}^{N2^{J+j}-1}\{ E[ (\hat{d}_{jk}-d_{jk})^{2}I(|\hat{d}_{jk}|>\delta_{j})] +E[ d_{jk}^{2}I(|\hat{d}_{jk}|\leq\delta_{j})] \}\nonumber
\\
&&{}+
\sum_{j=q+1}^{\infty}\sum_{k=-N+1}^{N2^{J+j}-1}d_{jk}^{2}=\Lambda_{1}+\Lambda_{2}+\Lambda_{3}+\Lambda_{4}.\nonumber
\end{eqnarray}
The proof then follows from Lemmas \ref{Lambda_1}--\ref{Lambda_6}, given below.
\end{pf*}

\begin{lemma}
\label{Lambda_1} Suppose that {the first derivatives of $g$ and
$\phi$ exist except for a finite number of points. {Moreover, assume
that $g^{\prime }$ and $\phi^{\prime}$ (where they exist) are
piecewise continuous and bounded.}} Then,
\[
\Lambda _{1}=\sum_{k=-N+1}^{N2^{J}-1}[ E(\hat{s}_{Jk})-s_{Jk}]
^{2}=\mathrm{O}(n^{-2}2^{2J}).
\]
\end{lemma}

\begin{pf}
For the expected value, we have
\begin{eqnarray*}
E(\hat{s}_{Jk})&=&E\Biggl( \frac{1}{n}\sum_{i=1}^{n}Y_{i}\phi_{Jk}(t_{i})%
\Biggr) =\frac{N^{1/2}2^{J/2}}{n}\sum_{i=1}^{n}g\biggl( \frac{i}{n}%
\biggr) \phi\biggl( N2^{J}\frac{i}{n}-k\biggr)
\\
&=&N^{1/2}2^{J/2}\sum_{i=1}^{n}\frac{1}{n}g\biggl( \frac{i}{n}\biggr) \phi\biggl(
N2^{J}\frac{i}{n}-k\biggr) .
\end{eqnarray*}
First, assume that $g$ and $\phi$ are continuously differentiable
and recall the rectangle rule
\begin{equation}\label{8}
\int_{a}^{b}f(t)\,\mathrm{d}t=\frac{b-a}{n}\sum_{i=0}^{n-1}f\biggl( a+i\frac{%
(b-a)}{n}\biggr)
+\mathrm{O}\Biggl( \sum_{i=0}^{n-1}\sup_{t\in I_{i}}|f^{\prime}(t)|\cdot\frac{(b-a)^{2}}{%
n^{2}}\Biggr)
\end{equation}
with $I_i=[a+i\frac{(b-a)}{n},a+(i+1)\frac{(b-a)}{n}]$. Noting that
the support of $\phi( N2^{J}t-k) $ (as a~function of $t$)
is $[kN^{-1}2^{-J},( kN^{-1}+1) 2^{-J}],$ we obtain
\[
E(\hat{s}_{Jk})=N^{1/2}2^{J/2}\sum_{i=i_{1}(k)}^{i_{2}(k)}\frac{1}{n}%
g\biggl( \frac{i}{n}\biggr) \phi\biggl( N2^{J}\frac{i}{n}-k\biggr)
\]
with
\[
i_1(k)= nkN^{-1}2^{-J}
\]
and
\[
i_2(k)=n( kN^{-1}+1) 2^{-J}.
\]
Thus, the number of non-zero terms in the sum is $n2^{-J}+1$. This,
together with the rectangle rule for $f(i/n)=g( i/n) \phi(
N2^{J}i/n-k) $ (and integration limits $a=0$, $b=1$), implies that
\[
E(\hat{s}_{Jk})=N^{1/2}2^{J/2}\int_{0}^{1}g(t)\phi(N2^{J}t-k)%
\,\mathrm{d}t+\mathrm{O}(n^{-1}2^{J/2})=s_{Jk}+\mathrm{O}(n^{-1}2^{J/2}).
\]
Note that, here, the factor $2^{J}$ from the derivative of $\phi(
N2^{J}t-k)$ is compensated by the fact that the number of non-zero
terms in the sum is proportional to $2^{-J}$.

Now, assume, more generally, that $g^{\prime }$ and $\phi^{\prime}$
exist except for a finite number of points and, where they exist,
that they are piecewise continuous and bounded. The result then follows
by a piecewise application of the rectangle rule.

In summary, we have
\[
E(\hat{s}_{Jk})-s_{Jk}=\mathrm{O}(n^{-1}2^{J/2}).
\]
This implies
that
\[
\Lambda _{1}=\sum_{k=-N+1}^{N2^{J}-1}[ E( \hat{s}_{Jk}) -s_{Jk}%
] ^{2}=\mathrm{O}\Biggl( \sum_{k=-N+1}^{N2^{J}-1}n^{-2}2^{J}\Biggr) =\mathrm{O}(n^{-2}2^{2J}),
\]
which completes the proof.
\end{pf}

\begin{lemma}
\label{Lambda_2}
Suppose that the first derivative of $\phi$ exists
on $[0,N]$ except for a finite number of points and, where
$\phi^{\prime}$ exists, it is piecewise continuous and bounded. Let
$J \ge 0 $ and $-N+1 \le k \le N2^{J}-1$. Then,
\[
E\{ [ \hat s_{Jk}-E( \hat s_{Jk}) ] ^{2}\} = C_{\phi}^{2}
N^{-1+\alpha}n^{-\alpha}2^{-J(1-\alpha)}+\mathrm{O}(n^{-1})
\]
and
\[
\Lambda_{2}=\sum_{k=-N+1}^{N2^{J}-1}E\{ [ \hat s_{Jk}-E( \hat s_{Jk}) ]
^{2}\} = C_{\phi}^{2} n^{-\alpha}N^{\alpha}2^{\alpha
J}+\mathrm{O}(n^{-1}2^{ J})+\mathrm{O}\bigl(n^{-\alpha}2^{-J(1-\alpha)}\bigr),
\]
where $C_{\phi}$ is the constant in (\ref{C_phi}).
\end{lemma}

\begin{pf}
First, assume that $\phi$ is continuously differentiable. Note that
$C_{\phi}$ is a positive finite constant (see Li and Xiao \cite{r39}).
We now consider the behavior of $E\{ [ \hat{s}_{Jk}-E( \hat{s}_{Jk})]^{2}\}$. We have
\begin{eqnarray*}
&&E\{ [ \hat{s}_{Jk}-E( \hat{s}_{Jk}) ]^{2}\}
\\
&&\quad=
E\Biggl\{ \Biggl[
\frac{1}{n}\sum_{i=1}^{n}\bigl( Y_{i}-E( Y_{i}) \bigr) \phi_{Jk}(t_{i})\Biggr] ^{2}\Biggr\}
\\
&&\quad=
E\Biggl[ \Biggl( \frac{N^{1/2}2^{J/2}}{n}\sum_{i=1}^{n}\xi_{i}\phi%
\biggl( N2^{J}\frac{i}{n}-k\biggr) \Biggr) ^{2}\Biggr]
\\
&&\quad=
Nn^{-2}2^{J}\sum_{i=1}^{n}\sum_{l=1}^{n}E(\xi_{i}\xi_{l})\phi%
\biggl( N2^{J}\frac{i}{n}-k\biggr) \phi\biggl( N2^{J}\frac{l}{n}-k\biggr)
\\
&&\quad=
Nn^{-2}2^{J}\sum_{i=nkN^{-1}2^{-J}}^{(1+kN^{-1})n2^{-J}}\sum%
_{l=nkN^{-1}2^{-J}}^{(1+kN^{-1})n2^{-J}}\gamma(l-i)\phi\biggl( N2^{J}%
\frac{i}{n}-k\biggr) \phi\biggl(
N2^{J}\frac{l}{n}-k\biggr)
\\
&&\quad=
Nn^{-2}2^{J}\mathop{\sum_{i,l=nkN^{-1}2^{-J}}}_{i\neq l}^{(1+kN^{-1})n2^{-J}}\gamma(l-i)\phi\biggl( N2^{J}\frac{i}{n}-k\biggr)
\phi\biggl( N2^{J}\frac{l}{n}-k\biggr)
\\
&&\qquad{}+
Nn^{-2}2^{J}\gamma(0)\sum_{i=nkN^{-1}2^{-J}}^{(1+kN^{-1})n2^{-J}}\phi^{2}\biggl( N2^{J}\frac{i}{n}-k\biggr) .
\end{eqnarray*}
Equation (\ref{8}) implies that
\begin{eqnarray*}
&&Nn^{-2}2^{J}\gamma(0)\sum_{i=nkN^{-1}2^{-J}}^{(1+kN^{-1})n2^{-J}}\phi^{2}\biggl( N2^{J}\frac{i}{n}-k\biggr)
\\
&&\quad=n^{-1}\gamma(0)\Biggl( \frac{N}{n2^{-J}}\sum_{i=nkN^{-1}2^{-J}}^{(1+kN^{-1})n2^{-J}}\phi^{2}\biggl( N2^{J}\frac{i}{n}-k\biggr) \Biggr)
\\
&&\quad=
n^{-1}\gamma(0)\int_{0}^{N}\phi^{2}(t)\,\mathrm{d}t+\mathrm{o}(n^{-1}).
\end{eqnarray*}
Due to (\ref{intphi2}), this is equal to
\[
n^{-1}\gamma(0)+\mathrm{o}(n^{-1})=\mathrm{O}(n^{-1}).
\]
Hence,
\begin{eqnarray*}
&&E\{ [ \hat{s}_{Jk}-E( \hat{s}_{Jk}) ] ^{2}\}
\\
&&\quad=
Nn^{-2}2^{J}\mathop{\sum_{i,l=nkN^{-1}2^{-J}}}_{i\neq l}^{(1+kN^{-1})n2^{-J}}\gamma(l-i)\phi\biggl( N2^{J}\frac{i}{n}-k\biggr) \phi\biggl(N2^{J}\frac{l}{n}-k\biggr) +\mathrm{O}(n^{-1}).
\end{eqnarray*}
Again using formula (\ref{5}), we obtain, by arguments analogous to those in, for example, Taqqu~\cite{r46},
\begin{eqnarray*}
&&E\{ [ \hat{s}_{Jk}-E( \hat{s}_{Jk}) ] ^{2}\}
\\
&&\quad\sim
C_{\gamma}Nn^{-2}2^{J}\mathop{\sum_{i,l=nkN^{-1}2^{-J}}}_{i\neq l}^{(1+kN^{-1})n2^{-J}}|l-i|^{-\alpha}\phi\biggl( N2^{J}\frac{i}{n}-k\biggr) \phi\biggl( N2^{J}\frac{l}{n}-k\biggr)
\\
&&\quad=
C_{\gamma}N^{\alpha}n^{-1-\alpha}2^{\alpha J}\sum_{i=nkN^{-1}2^{-J}}^{(1+kN^{-1})n2^{-J}}\phi\biggl( N2^{J}\frac{i}{n}-k\biggr)\frac{N2^{J}}{n}
\\
&&\quad\hphantom{=C_{\gamma}N^{\alpha}n^{-1-\alpha}2^{\alpha J}\sum_{i=nkN^{-1}2^{-J}}^{(1+kN^{-1})n2^{-J}}}{}\times
\mathop{\sum_{l=nkN^{-1}2^{-J}}}_{l\neq i}^{(1+kN^{-1})n2^{-J}}\biggl\vert N2^{J}\frac{l}{n}-N2^{J}\frac{i}{n}\biggr\vert ^{-\alpha}\phi\biggl( N2^{J}\frac{l}{n}-k\biggr) .
\end{eqnarray*}
The function $f(x)=\vert x-( N2^{J}\frac{i}{n}-k)
\vert ^{-\alpha}\phi( x) $ is differentiable on $[0,N2^{J}\frac{%
i-1}{n}-k]\cup[N2^{J}\frac{i+1}{n}-k,N]$ for all fixed $i$ and
$n$. Therefore, the rectangle rule implies that
\begin{eqnarray*}
&&
\frac{N}{n2^{-J}}\mathop{\sum_{l=nkN^{-1}2^{-J}}}_{l\neq i}^{(1+kN^{-1})n2^{-J}}\biggl\vert N2^{J}\frac{l}{n}-N2^{J}\frac{i}{n}\biggr\vert ^{-\alpha}\phi\biggl( N2^{J}\frac{l}{n}-k\biggr)
\\
&&\quad=
\frac{N}{n2^{-J}}\sum_{l=nkN^{-1}2^{-J}}^{i-1}\biggl\vert\biggl(N2^{J}\frac{l}{n}-k\biggr) -\biggl( N2^{J}\frac{i}{n}-k\biggr) \biggr\vert ^{-\alpha}\phi\biggl(N2^{J}\frac{l}{n}-k\biggr)
\\
&&\qquad{}+
\frac{N}{n2^{-J}}\sum_{l=i+1}^{(1+kN^{-1})n2^{-J}}\biggl\vert \biggl(N2^{J}\frac{l}{n}-k\biggr) -\biggl( N2^{J}\frac{i}{n}-k\biggr)\biggr\vert^{-\alpha}\phi\biggl(N2^{J}\frac{l}{n}-k\biggr)
\\
&&\quad=
\int_{0}^{N2^{J}((i-1)/n)-k}\biggl\vert x-\biggl( N2^{J}\frac{i}{n}-k\biggr) \biggr\vert ^{-\alpha}\phi(x)\,\mathrm{d}x
\\
&&\qquad{}+
\int_{N2^{J}((i+1)/n)-k}^{N}\biggl\vert x-\biggl( N2^{J}\frac{i}{n}-k\biggr) \biggr\vert ^{-\alpha}\phi(x)\,\mathrm{d}x+K_{1,n}+K_{2,n},
\end{eqnarray*}
where
\[
K_{1,n}=\mathrm{O}\Biggl( \biggl( \frac{N}{n2^{-J}}\biggr)^{2}
\sum_{l=nkN^{-1}2^{-J}}^{i-2}\sup_{x\in I_{l}(k)}\biggl\vert \frac{%
\mathrm{d}}{\mathrm{d}x}\biggl( \biggl\vert x-\biggl( N2^{J}\frac{i}{n}-k\biggr) \biggr\vert ^{-\alpha}\phi(
x) \biggr) \biggr\vert \Biggr)
\]
with $I_{l}(k)=[N2^{J}l/n-k,N2^{J}(l+1)/n-k]$ and
\begin{eqnarray*}
K_{2,n}&=&\mathrm{O}\Biggl( \biggl( \frac{N}{n2^{-J}}\biggr)
^{2}
\\
&&\hphantom{\mathrm{O}\Biggl(}{}\times\sum_{l=i+1}^{(1+kN^{-1})n2^{-J}-1}\sup_{x\in I_{l}(k)}\biggl\vert
\frac{\mathrm{d}}{\mathrm{d}x}\biggl( \biggl\vert x-\biggl( N2^{J}\frac {i}{n}-k\biggr) \biggr\vert ^{-\alpha}\phi( x)
\biggr) \biggr\vert \Biggr) .
\end{eqnarray*}
Now,
\begin{eqnarray*}
&&
\Biggl\vert \biggl( \frac{N}{n2^{-J}}\biggr)
^{2}\sum_{l=nkN^{-1}2^{-J}}^{i-2}\sup_{x\in I_{l}(k)}\biggl\vert \frac{%
\mathrm{d}}{\mathrm{d}x}\biggl( \biggl\vert x-\biggl( N2^{J}\frac{i}{n}-k\biggr) \biggr\vert ^{-\alpha}\phi(
x) \biggr) \biggr\vert \Biggr\vert
\\
&&\quad\leq
\alpha N^2\max_{x\in\lbrack0,1]}\phi( x) \cdot\Biggl\vert
n^{-2}2^{2J}\sum_{l=nkN^{-1}2^{-J}}^{i-2}\biggl( \biggl( N2^{J}\frac {i%
}{n}-k\biggr) -\biggl( N2^{J}\frac{l+1}{n}-k\biggr) \biggr) ^{-1-\alpha }\Biggr\vert
\\
&&\qquad{}+
N^2\max_{x\in\lbrack0,1]}\phi^{\prime}( x) \cdot\Biggl\vert
n^{-2}2^{2J}\sum_{l=nkN^{-1}2^{-J}}^{i-2}\biggl( \biggl( N2^{J}\frac {i%
}{n}-k\biggr) -\biggl( N2^{J}\frac{l+1}{n}-k\biggr) \biggr) ^{-\alpha }\Biggr\vert
\\
&&\quad\leq
C_{1}n^{-(1-\alpha)}2^{J(1-\alpha)}\sum%
_{j=1}^{i-1-nkN^{-1}2^{-J}}j^{-1-\alpha}
+
C_{2}n^{-(2-\alpha)}2^{J(2-\alpha)}\sum%
_{j=1}^{i-1-nkN^{-1}2^{-J}}j^{-\alpha}
\\
&&\quad\leq
C_{1}n^{-(1-\alpha)}2^{J(1-\alpha)}\sum_{j=1}^{\infty}j^{-1-%
\alpha}+C_{2}n^{-(2-\alpha)}2^{J(2-\alpha)}\sum_{j=1}^{n2^{-J}}j^{-%
\alpha}
\\
&&\quad\leq
C_{1}n^{-(1-\alpha)}2^{J(1-\alpha)}\sum_{j=1}^{\infty}j^{-1-%
\alpha}+C_{2}^{\ast}n^{-1}2^{J}.
\end{eqnarray*}
Thus,
\[
K_{1,n}=\mathrm{O}\bigl(n^{-(1-\alpha )}2^{J(1-\alpha)}\bigr).
\]
By analogous arguments, we obtain%
\[
K_{2,n}=\mathrm{O}\bigl(n^{-(1-\alpha )}2^{J(1-\alpha)}\bigr).
\]
This implies that
\begin{eqnarray*}
&&E\{ [ \hat{s}_{Jk}-E( \hat{s}_{Jk})] ^{2}\}
\\[-2pt]
&&\quad=C_{\gamma}N^{\alpha}n^{-1-\alpha}2^{\alpha J}
\\[-2pt]
&&\qquad{}\times
\sum_{i=nkN^{-1}2^{-J}}^{(1+kN^{-1})n2^{-J}}\phi\biggl( N2^{J}\frac{i}{n}-k\biggr)
\biggl( \int_{0}^{N2^{J}((i-1)/n)-k}\biggl\vert x-\biggl( N2^{J}\frac{i}{n}-k\biggr) \biggr\vert ^{-\alpha}\phi( x)\,\mathrm{d}x
\\[-2pt]
&&\qquad\hphantom{{}\times\sum_{i=nkN^{-1}2^{-J}}^{(1+kN^{-1})n2^{-J}}\phi\biggl( N2^{J}\frac{i}{n}-k\biggr)\biggl(}{}+
\int_{N2^{J}((i+1)/n)-k}^{N}\biggl\vert x-\biggl( N2^{J}\frac {i}{n}-k\biggr) \biggr\vert ^{-\alpha}\phi( x)\,\mathrm{d}x
\\[-2pt]
&&\qquad\hphantom{{}\times\sum_{i=nkN^{-1}2^{-J}}^{(1+kN^{-1})n2^{-J}}\phi\biggl( N2^{J}\frac{i}{n}-k\biggr)\biggl(}{}+
\mathrm{O}\bigl(n^{-(1-\alpha)}2^{J(1-\alpha)}\bigr)\biggr) +\mathrm{O}(n^{-1})
\\[-2pt]
&&\quad=
C_{\gamma}N^{\alpha}n^{-1-\alpha}2^{\alpha J}\sum_{i=nkN^{-1}2^{-J}}^{(1+kN^{-1})n2^{-J}}
\phi\biggl( N2^{J}\frac{i}{n}-k\biggr)
\\[-2pt]
&&\quad\hphantom{=C_{\gamma}N^{\alpha}n^{-1-\alpha}2^{\alpha J}\sum_{i=nkN^{-1}2^{-J}}^{(1+kN^{-1})n2^{-J}}}{}\times
\int_{0}^{N2^{J}((i-1)/n)-k}\biggl\vert x-\biggl( N2^{J}\frac{i}{n}-k\biggr) \biggr\vert ^{-\alpha}\phi( x)\,\mathrm{d}x
\\[-2pt]
&&\qquad{}+
C_{\gamma}N^{\alpha}n^{-1-\alpha}2^{\alpha J}\sum_{i=nkN^{-1}2^{-J}}^{(1+kN^{-1})n2^{-J}}
\phi\biggl( N2^{J}\frac{i}{n}-k\biggr)
\\[-2pt]
&&\qquad\hphantom{{}+C_{\gamma}N^{\alpha}n^{-1-\alpha}2^{\alpha J}\sum_{i=nkN^{-1}2^{-J}}^{(1+kN^{-1})n2^{-J}}}{}\times
\int _{N2^{J}((i+1)/n)-k}^{N}\biggl\vert x-\biggl( N2^{J}\frac{i}{n}-k\biggr) \biggr\vert ^{-\alpha}\phi( x)\,\mathrm{d}x
\\[-2pt]
&&\qquad{}+
\mathrm{O}(n^{-1})=A_{1}+A_{2}+\mathrm{O}(n^{-1}).
\end{eqnarray*}
Again using (\ref{8}), we obtain, by arguments analogous to those used above,
\[
A_{1}=C_{\gamma}N^{-1+\alpha}n^{-\alpha}2^{-J(1-\alpha)}
\int_{0}^{N}\!\!\!\int_{0}^{y-N2^{J}(1/n)}|x-y|^{-\alpha}\phi(x)\phi(y)\,\mathrm{d}x\,\mathrm{d}y+\mathrm{O}(n^{-1})
\]
and
\[
A_{2}=C_{\gamma}N^{-1+\alpha}n^{-\alpha}2^{-J(1-\alpha)}\int_{0}^{N}\!\!\!\int_{y+N2^{J}(1/n)}^{N}|x-y|^{-\alpha}\phi(x)\phi
(y)\,\mathrm{d}x\,\mathrm{d}y+\mathrm{O}(n^{-1}).
\]
Noting that
\[
\int_{y-N2^{J}(1/n)}^{y+N2^{J}(1/n)}\!\vert x\,{-}\,y\vert
^{-\alpha}\phi( x)
\phi(y)\,\mathrm{d}x\,{\le}\,2\max_{x \in [0,N]}(\phi^2( x))\,{\cdot}\,%
\int_{0}^{N2^{J}n^{-1}}\!z^{-\alpha}\,\mathrm{d}y\,{=}\,\mathrm{O}\bigl(n^{-(1-\alpha)}2^{J(1-%
\alpha)}\bigr),
\]
we obtain
\[
\int_{0}^{N}\!\!\!\int_{y-N2^{J}(1/n)}^{y+N2^{J}(1/n)}\vert
x-y\vert ^{-\alpha}\phi( x) \phi(y)\,\mathrm{d}x\,\mathrm{d}y=\mathrm{O}\bigl(n^{-(1-\alpha
)}2^{J(1-\alpha)}\bigr)\vadjust{\goodbreak}
\]
and
\begin{eqnarray*}
E\{ [ \hat{s}_{Jk}-E( \hat{s}_{Jk}) ] ^{2}\}
&=&C_{\gamma}N^{-1+\alpha}n^{-\alpha}2^{-J(1-\alpha)}\int_{0}^{N}\!\!\!\int%
_{0}^{N}|x-y|^{-\alpha}\phi(x)\phi (y)\,\mathrm{d}x\,\mathrm{d}y+\mathrm{O}(n^{-1})
\\
&=&C_{\phi}^{2}N^{-1+\alpha}n^{-\alpha}2^{-J(1-\alpha)}+\mathrm{O}(n^{-1}),
\end{eqnarray*}
where $C_{\phi}$ is the constant in (\ref{C_phi}). Hence,
\begin{eqnarray*}
\Lambda_{2}&=&\sum_{k=-N+1}^{N2^{J}-1}E\{ [ \hat{s}_{Jk}-E( \hat{s}_{Jk})
] ^{2}\} =\sum _{k=-N+1}^{N2^{J}-1}\bigl(
C_{\phi}^{2}N^{-1+\alpha}n^{-\alpha }2^{-J(1-\alpha)}+\mathrm{O}(n^{-1})\bigr)
\\
&=&C_{\phi}^{2}n^{-\alpha}N^{\alpha}2^{\alpha
J}+\mathrm{O}(n^{-1}2^{J})+\mathrm{O}\bigl(n^{-\alpha }2^{-J(1-\alpha)}\bigr).
\end{eqnarray*}

In the general case where $\phi^{\prime}$ exists except for a finite
number of points and, where it exists, it is piecewise continuous
and bounded, the result follows by a piecewise  application of the rectangle rule.
\end{pf}

\begin{lemma}
\label{Lambda_3}
Suppose that the first derivative of $\psi$ exists
on $[0,N]$ except for a finite number of points and, where
$\psi^{\prime}$ exists, it is piecewise continuous and bounded. Let
$J \ge0 $, $j\ge0$ and $-N+1 \le k \le N2^{J+j}-1$. Then
\begin{eqnarray*}
\sigma_{j}^{2}&=&E\{ [ \hat d_{jk}-E(\hat d_{jk})]^{2}\}
\\
&=&C_{\psi}^{2}N^{-1+\alpha}n^{-\alpha}2^{-(J+j)(1-\alpha)}+\mathrm{O}(n^{-1}),
\end{eqnarray*}
where $C_{\psi}$ is the constant in (\ref{C_psi}).
\end{lemma}

\begin{pf}
Noting that
\[
E\{ [ \hat{d}_{jk}-E(\hat{d}_{jk})] ^{2}\} =E\Biggl\{ %
\Biggl[ \frac{N^{1/2}2^{(J+j)/2}}{n}\sum_{i=1}^{n}\xi _{i}\psi\biggl(
N2^{J+j}\frac{i}{n}-k\biggr) \Biggr] ^{2}\Biggr\},
\]
the proof is analogous to the proof of Lemma \ref{Lambda_2}, {with the difference being that $%
\psi$ is used instead of $\phi$  and $J$ is replaced by $J+j$.}
\end{pf}

\begin{lemma}
\label{Lambda_4}
Suppose that the first $r$ derivatives of $g$
exist and are continuous on $[0,1]$. Then, for all $j\geq0$ and
$0\leq k\leq N2^{J+j}-1$,
\begin{eqnarray}\label{djk}
d_{jk}&=&\frac{\nu_{r}}{r!}g^{(r)}\bigl(kN^{-1}2^{-(J+j)}\bigr)N^{-(2r+1)/2}2^{-((2r+1)/2)(J+j)}\nonumber
\\[-8pt]\\[-8pt]
&&{}+\mathrm{o}\bigl( 2^{-((2r+1)/2)(J+j)}\bigr),\nonumber
\end{eqnarray}
where $\nu_r$ is the $r$th moment of $\psi$ (see (\ref{10})).
Together with the assumptions of Lemma \ref{Lambda_3}, this yields
that
\[
E(\hat{d}_{jk})-d_{jk}=\mathrm{O}\bigl(n^{-1}2^{(J+j)/2}\bigr).
\]
\end{lemma}

\begin{pf}
Note that
\begin{eqnarray*}
d_{jk}&=&N^{1/2}2^{(J+j)/2}\int_{0}^{1}g( t)
\psi(N2^{J+j}t-k)\,\mathrm{d}t
\\
&=&N^{1/2}2^{(J+j)/2}\int_{kN^{-1}2^{-(J+j)}}^{(1+kN^{-1})2^{-(J+j)}}g\bigl(
N^{-1}2^{-(J+j)}[ N2^{J+j}t-k+k] \bigr) \psi(N2^{J+j}t-k)\,\mathrm{d}t
\\
&=&N^{-1/2}2^{-(J+j)/2}\int_{0}^{N}g\bigl(N^{-1}2^{-(J+j)}(y+k)\bigr)\psi(y)\,\mathrm{d}y.
\end{eqnarray*}
Since $g$ is $r$-times continuously differentiable, {the local
Taylor} expansion {(see, e.g., Zorich~\cite{r52}, pages 225--226)} of $g$ yields
\begin{eqnarray*}
d_{jk}
&=&
N^{-1/2}2^{-(J+j)/2}\int_{0}^{N}\psi(y)\biggl[ g\bigl(kN^{-1}2^{-(J+j)}\bigr)+N^{-1}2^{-(J+j)}g^{\prime}\bigl(kN^{-1}2^{-(J+j)}\bigr)y
\\
&&\hphantom{N^{-1/2}2^{-(J+j)/2}\int_{0}^{N}\psi(y)\biggl[}{}+
\cdots+\frac{N^{-r}2^{-r(J+j)}}{r!}g^{(r)}\bigl(kN^{-1}2^{-(J+j)}\bigr)y^{r}\biggr]\,\mathrm{d}y
\\
&&{}+\mathrm{o}\bigl( 2^{-((2r+1)/2)(J+j)}\bigr).
\end{eqnarray*}
The moment conditions (\ref{9}) and (\ref{10}) then imply
that
\begin{eqnarray*}
d_{jk}&\!=\!&\frac{1}{r!}g^{(r)}\bigl(kN^{-1}2^{-(J+j)}\bigr)N^{-(2r+1)/2}2^{-((2r+1)/2)(J+j)}\!\!\int_{0}^{N}\!y^{r}\psi(y)\,\mathrm{d}y
\,{+}\,\mathrm{o}\bigl( 2^{-((2r+1)/2)(J+j)}\bigr)
\\
&\!=\!&\frac{\nu_{r}}{r!}g^{(r)}\bigl(kN^{-1}2^{-(J+j)}\bigr)N^{-(2r+1)/2}2^{-((2r+1)/2)(J+j)}\,{+}\,\mathrm{o}\bigl( 2^{-((2r+1)/2)(J+j)}\bigr) .
\end{eqnarray*}
For $E(\hat{d}_{jk}),$ we have
\begin{eqnarray*}
E(\hat{d}_{jk})&=&\frac{1}{n}\sum_{i=1}^{n}E[ Y_{i}\psi_{jk}(t_{i})]
\\
&=&N^{1/2}2^{(J+j)/2}\sum_{i=1}^{n}n^{-1}g\biggl( \frac{i}{n}\biggr) \psi\biggl(
N2^{J+j}\frac{i}{n}-k\biggr) .
\end{eqnarray*}
Again using {the same arguments as in Lemma \ref{Lambda_1} for
$E(\hat{s}_{Jk})$}, we obtain
that
\[
E(\hat{d}_{jk})=d_{jk}+\mathrm{O}\bigl(n^{-1}2^{(J+j)/2}\bigr).
\]
\upqed\end{pf}

\begin{lemma}
\label{Lambda_5}  Under the assumptions of Lemma \ref{Lambda_4},
\[
\Lambda_{4}=\frac{1}{(r!)^{2}}\frac{1}{2^{2r}-1}N^{-2r}2^{-2r(J+q)}\int%
_{0}^{1}\nu_{r}^{2}\bigl(g^{(r)}(t)\bigr)^2\,\mathrm{d}t+\mathrm{o}\bigl( 2^{-2r(J+q)}\bigr) .
\]
\end{lemma}

\begin{pf}
Using (\ref{MISE}), %(24)
we have
\begin{eqnarray*}
\Lambda_{4} & =&\sum_{j=q+1}^{\infty}\sum_{k=-N+1}^{N2^{J+j}-1}d_{jk}^{2} \\
& =&\sum_{j=q+1}^{\infty}\sum_{k=-N+1}^{N2^{J+j}-1}\biggl[ \frac{\nu_{r}}{r!}%
g^{(r)}\bigl(kN^{-1}2^{-(J+j)}\bigr)\biggr] ^{2}N^{-(2r+1)}2^{-(2r+1)(J+j)}+\mathrm{o}\bigl(
2^{-(2r+1)(J+j)}\bigr) .
\end{eqnarray*}
Note that the continuity of $g^{(r)}$ implies convergence of the Riemann
sum. Hence, $\Lambda_4$ is
equal to
\begin{eqnarray*}
&&
\frac{1}{(r!)^{2}}\sum_{j=q+1}^{\infty}N^{-2r}2^{-2r(J+j)}\biggl\{
\int_{0}^{1}\nu_{r}^{2}\bigl(g^{(r)}(t)\bigr)^2\,\mathrm{d}t+\mathrm{o}(1)\biggr\} +\mathrm{o}\bigl( 2^{-2r(J+q)}\bigr)
 \\
&&\quad=
\frac{1}{(r!)^{2}}\frac{1}{2^{2r}-1}N^{-2r}2^{-2r(J+q)}\int%
_{0}^{1}\nu_{r}^{2}\bigl(g^{(r)}(t)\bigr)^2\,\mathrm{d}t+\mathrm{o}\bigl( 2^{-2r(J+q)}\bigr).
\end{eqnarray*}
\upqed\end{pf}

\begin{lemma}
\label{Lambda_6}
Let
\begin{equation}\label{Hatq}
\hat{q}=\log_{2}n^{\alpha/(2r+\alpha)}+\frac{1}{2r+\alpha}\log
_{2}\biggl( \frac{\nu_{r}^{2}}{(r!)^{2}C_{\psi}^{2}N^{2r+\alpha}}\max
_{t\in[ 0,1]}\bigl[ g^{(r)}(t)\bigr] ^{2}\biggr) -J+1
\end{equation}
and
\[
\lambda_{jk}=E[ (\hat{d}_{jk}-d_{jk})^{2}I(|\hat{d}_{jk}|>\delta _{j})%
] +E[ d_{jk}^{2}I(|\hat{d}_{jk}|\leq\delta_{j})] .
\]
Under the assumptions of Lemmas \ref{Lambda_3} and \ref{Lambda_4}, the
following then holds: If $q>\hat{q}$, then for all $j$ with  $\hat{q}<
j < \frac{2+\alpha}{4r+2+\alpha}\log_2n-J$, we have
\[
\min_{\delta_{j}}\lambda_{jk}=d_{jk}^{2}+\mathrm{O}\bigl(
2^{\alpha(J+j)/2}n^{-(1+\alpha/2)}\bigr).
\]
\end{lemma}

\begin{pf}
Defining
$S_0=2^{-(2r+\alpha)}C_{\psi}^{2}N^{-1+\alpha}n^{-\alpha}2^{-(J+j)(1-\alpha)}$
and taking into account Lemma \ref{Lambda_3}, we have
$S_{1,j}=2^{-(2r+\alpha)}\sigma_j^2S_{0}^{-1}=1+r_{1,j}$ with
$|r_{1,j}|\le r_1=\mathrm{o}(1)$ for all $j\ge \hat q$. Moreover, Lemma
\ref{Lambda_4} implies that
\begin{eqnarray*}
S_{2,jk}&=&d_{jk}^{2}S_0^{-1}=\frac{\nu_{r}^{2}}{(r!)^{2}C_{\psi}^{2}N^{2r+\alpha}}\bigl[
g^{(r)}\bigl(kN^{-1}2^{-(J+j)}\bigr)\bigr]^22^{-(2r+\alpha)(J+j-1)}n^{\alpha}+r_{2,jk}
\\
&\le& \frac{\nu_{r}^{2}}{(r!)^{2}C_{\psi}^{2}N^{2r+\alpha}}\max _{t\in
[0,1]}\bigl[ g^{(r)}(t)\bigr] ^{2}n^{\alpha}\max_{j>\hat q}\bigl\{
2^{-(2r+\alpha)(J+j-1)}\bigr\}+r_2
\end{eqnarray*}
with $r_2=\mathrm{o}(1)$ independent of $j$ and
$k$.
Using (\ref{Hatq}), we obtain
\[
S_{2,jk}\le 2^{-(2r+\alpha)}+r_2 \le 1+r_1=S_{1,j}
\]
for $j>\hat q$ and $n$ large enough, which implies
that%
\begin{equation}\label{Relation}
\sigma_{j}\geq2^{r+\alpha/2}\max_{k}|d_{jk}|.
\end{equation}
The mean squared error $\lambda_{jk}$ can be written as
\begin{eqnarray*}
\lambda_{jk}
& =&
E[ (\hat{d}_{jk}-d_{jk})^{2}I(|\hat{d}_{jk}|>\delta_{j})%
] +E[ d_{jk}^{2}I(|\hat{d}_{jk}|\leq\delta _{j})]
\\
& =&
\frac{1}{\sqrt{2\uppi}\sigma_{j}}\int_{|t|>\delta_{j}}\!(t-d_{jk})^{2}\mathrm{e}^{-(1/(2\sigma_{j}^{2}))(t-E(\hat{d}_{jk}))^{2}}\,\mathrm{d}t+d_{jk}^{2}\frac{1}{%
\sqrt{2\uppi}\sigma_{j}}\int_{|t|<\delta_{j}}\!\mathrm{e}^{-(t-E(\hat {d}_{jk}))^{2}%
/(2\sigma_{j}^{2})}\,\mathrm{d}t
\\
& =&
A_{1}+A_{2}.
\end{eqnarray*}
We approximate $A_{1}$ and $A_{2}$ separately. Taylor expansion of
$A_{1}$ with respect to $E(\hat{d}_{jk})$ in the neighborhood of
$d_{jk}$ yields
\begin{eqnarray*}
A_{1}
&=&
\frac{1}{\sqrt{2\uppi}\sigma_{j}}\int_{|t|>\delta_{j}}(t-d_{jk})^{2}\mathrm{e}^{-(1/(2\sigma_{j}^{2}))(t-E(\hat{d}_{jk}))^{2}}\,\mathrm{d}t
\\
&=&
\frac{1}{\sqrt{2\uppi}\sigma_{j}}\int_{|t|>\delta_{j}}(t-d_{jk})^{2}\mathrm{e}^{-(1/(2\sigma_{j}^{2}))(t-d_{jk})^{2}}\,\mathrm{d}t
\\
&&{}+
\frac{E(\hat{d}_{jk})-d_{jk}}{\sqrt{2\uppi}\sigma_{j}}\int_{|t|>\delta_{j}}%
\frac{(t-d_{jk})^{3}}{\sigma_{j}^{2}}\mathrm{e}^{-(1/(2\sigma_{j}^{2}))(t-d_{jk})^{2}}\,\mathrm{d}t
\\
&&{}+
\mathrm{O}\biggl( \frac{[ E(\hat{d}_{jk})-d_{jk}] ^{2}}{\sigma_{j}}%
\int_{|t|>\delta_{j}}\biggl(\frac{(t-d_{jk})^{4}}{\sigma_{j}^{4}}-\frac{(t-d_{jk})^{2}}{\sigma_{j}^{2}}\biggr)\mathrm{e}^{-(1/(2\sigma_{j}^{2}))(t-d_{jk})^{2}}\,\mathrm{d}t\biggr) .
\end{eqnarray*}
If $d_{jk}\neq 0$, then Lemmas \ref{Lambda_3} and \ref{Lambda_4} imply that
\begin{eqnarray}\label{RemTerm2}
&&\frac{E(\hat{d}_{jk})-d_{jk}}{\sigma_{j}}\int_{|t|>\delta_{j}}\frac {%
(t-d_{jk})^{3}}{\sigma_{j}^{2}}\mathrm{e}^{-(1/(2\sigma_{j}^{2}))(t-d_{jk})^{2}}\,\mathrm{d}t\nonumber
\\
&&\quad=\sigma_{j}[ E(\hat{d}_{jk})-d_{jk}] \int_{|t|>\delta_{j}/%
\sigma_{j}}\biggl( t-\frac{d_{jk}}{\sigma_{j}}\biggr) ^{3}\mathrm{e}^{-1/2
( t-d_{jk}/\sigma_{j}) ^{2}}\,\mathrm{d}t
\\
&&\quad=\mathrm{O}\bigl( n^{-\alpha/2}2^{-(J+j)(1-\alpha)/2}\cdot n^{-1}2^{(J+j)/2}\bigr) =\mathrm{O}\bigl(
2^{\alpha(J+j)/2}n^{-(1+\alpha/2)}\bigr) .\nonumber
\end{eqnarray}
If $d_{jk} = 0$, then
\[
\frac{E(\hat{d}_{jk})-d_{jk}}{\sigma_{j}}\int_{|t|>\delta_{j}}\frac {%
(t-d_{jk})^{3}}{\sigma_{j}^{2}}\mathrm{e}^{-(1/(2\sigma_{j}^{2}))(t-d_{jk})^{2}}\,\mathrm{d}t=0
\]
and
\begin{eqnarray}\label{RemTerm1}
&&\frac{[ E(\hat{d}_{jk})-d_{jk}] ^{2}}{\sigma_{j}}\int
_{|t|>\delta_{j}}\biggl(\frac{(t-d_{jk})^{4}}{\sigma_{j}^{4}}-\frac{(t-d_{jk})^{2}}{\sigma_{j}^{2}}\biggr)
\mathrm{e}^{-(1/(2\sigma_{j}^{2}))(t-d_{jk})^{2}}\,\mathrm{d}t\nonumber
\\
&&\quad=[ E(\hat{d}_{jk})-d_{jk}] ^{2}\int_{|t|>\delta_{j}/%
\sigma_{j}}\biggl(\biggl( t-\frac{d_{jk}}{\sigma_{j}}\biggr) ^{4}-\biggl( t-\frac{d_{jk}}{\sigma_{j}}\biggr) ^{2}\biggr)
\mathrm{e}^{-(1/2)( t-d_{jk}/\sigma_{j})^{2}}\,\mathrm{d}t
\\[-2pt]
&&\quad=\mathrm{O}( n^{-2}2^{J+j}) .\nonumber
\end{eqnarray}

The condition $j+J< \frac{2+\alpha}{4r+2+\alpha}\log_2n$ implies that ${n}%
^{-2}2^{J+j}=\mathrm{o}(2^{(\alpha/2)(J+j)}n^{-(1+\alpha/2)})$ so that
\begin{equation}\label{A1}
A_{1}=\frac{1}{\sqrt{2\uppi}\sigma_{j}}\int_{|t|>\delta_{j}}(t-d_{jk})^{2}\mathrm{e}^{(-1/(2\sigma_{j}^{2}))(t-d_{jk})^{2}}\,\mathrm{d}t+\mathrm{O}\bigl( 2^{\alpha
(J+j)/2}n^{-(1+\alpha/2)}\bigr).
\end{equation}

By analogous arguments, we have, for $d_{jk}\neq0$,
\begin{eqnarray*}
A_{2}&=&d_{jk}^{2}\frac{1}{\sqrt{2\uppi}\sigma_{j}}\int_{|t|<\delta_{j}}
\mathrm{e}^{-(t-E(\hat{d}_{jk}))^{2}/(2\sigma_{j}^{2})}\,\mathrm{d}t
\\[-2pt]
&=&
d_{jk}^{2}\frac{1}{\sqrt{2\uppi}\sigma_{j}}\int_{|t|<\delta_{j}}\mathrm{e}^{-(t-d_{jk})^{2}/(2\sigma_{j}^{2})}\,\mathrm{d}t
\\[-2pt]
&&{}+\mathrm{O}\biggl( \frac{d_{jk}^{2}[ E(\hat{d}_{jk})-d_{jk}] }{\sigma_{j}}\int_{|t|<\delta_{j}}\frac{(t-d_{jk})}{%
\sigma_{j}^{2}}\mathrm{e}^{-(t-d_{jk})^{2}/(2\sigma_{j}^{2})}\,\mathrm{d}t\biggr)
\end{eqnarray*}
with
\begin{eqnarray}\label{RemTerm3}
&&\frac{d_{jk}^{2}[ E(\hat{d}_{jk})-d_{jk}] }{\sigma_{j}}%
\int_{|t|<\delta_{j}}\frac{(t-d_{jk})}{\sigma_{j}^{2}}\mathrm{e}^{-(t-d_{jk})^{2}/(2\sigma_{j}^{2})}\,\mathrm{d}t\nonumber
\\[-2pt]
&&\quad=\frac{d_{jk}^{2}[ E(\hat{d}_{jk})-d_{jk}] }{\sigma_{j}}%
\int_{|t|<\delta_{j}/\sigma_{j}}\biggl( t-\frac{d_{jk}}{\sigma_{j}}\biggr)
\mathrm{e}^{-(1/2)( t-d_{jk}/\sigma_{j}) ^{2}}\,\mathrm{d}t\nonumber
\\[-9pt]\\[-9pt]
&&\quad=\mathrm{O}\bigl( 2^{-(2r+1)(J+j)}\cdot n^{-1}2^{(J+j)/2}\cdot n^{\alpha
/2}2^{(J+j)(1-\alpha)/2}\bigr)\nonumber
\\[-2pt]
&&\quad=\mathrm{O}\bigl( n^{-(1-\alpha/2)}2^{-(2r+\alpha/2)(J+j)}\bigr) .\nonumber
\end{eqnarray}

For $d_{jk}=0$, we have
\[
A_{2}=d_{jk}^{2}\frac{1}{\sqrt{2\uppi}\sigma_{j}}\int_{|t|<\delta_{j}}\mathrm{e}^{-(t-E(\hat{d}_{jk}))^{2}/(2\sigma_{j}^{2})}\,\mathrm{d}t=0.
\]
In summary, we have derived the approximation,
\begin{eqnarray}\label{RemTerm4}
\lambda_{jk}
&=&
A_{1}+A_{2}\nonumber
\\[-2pt]
&=&
\sigma_{j}^{2}\frac{1}{\sqrt{2\uppi}}\int_{|t|>\delta_{j}/\sigma_{j}}\biggl( t-%
\frac{d_{jk}}{\sigma_{j}}\biggr) ^{2}\mathrm{e}^{-(1/2)(
t-d_{jk}/\sigma_{j}) ^{2}}\,\mathrm{d}t\nonumber
\\[-9pt]\\[-9pt]
&&{}+
d_{jk}^{2}\frac{1}{\sqrt{2\uppi}}\int
_{|t|<\delta_{j}/\sigma_{j}}\mathrm{e}^{-(1/2)( t-d_{jk}/\sigma_{j}) ^{2}}\,\mathrm{d}t\nonumber
\\[-2pt]
&&{}+
\mathrm{O}\bigl( 2^{\alpha(J+j)/2}n^{-(1+\alpha/2)}\bigr) +\mathrm{O}\bigl( n^{-(1-\alpha
/2)}2^{-(2r+\alpha/2)(J+j)}\bigr)\nonumber
\end{eqnarray}
with uniformly bounded error terms (see (\ref{RemTerm2}),
 (\ref{RemTerm1}) and (\ref{RemTerm3})). It is then sufficient to show
that for all $k$ and  all $j$ with $\hat{q}< j <
\frac{2+\alpha}{4r+2+\alpha}\log_2n-J,$ we have
\[
\min_{\delta_{j}}\hat \lambda_{jk}=d_{jk}^{2},
\]
where
\begin{eqnarray}\label{hatlambda}
\hat{\lambda}_{jk}
&=&
\sigma_{j}^{2}\frac{1}{\sqrt{2\uppi}}\int_{|t|>\delta_{j}/\sigma_{j}}\biggl( t-\frac{d_{jk}}{\sigma_{j}}\biggr) ^{2}
\mathrm{e}^{-(1/2)( t-d_{jk}/\sigma_{j}) ^{2}}\,\mathrm{d}t\nonumber
\\[-8pt]\\[-8pt]
&&{}+
d_{jk}^{2}\frac{1}{\sqrt{2\uppi}}\int_{|t|<\delta_{j}/\sigma_{j}}\mathrm{e}^{-(1/2)( t-d_{jk}/\sigma_{j}) ^{2}}\,\mathrm{d}t.\nonumber
\end{eqnarray}
In the following, we distinguish two cases: $\delta_{j}\leq\sigma_{j}$ and $%
\delta_{j}>\sigma_{j}$.

At first, let $\delta_{j}\leq\sigma_{j}$. Recall that  $\sigma_{j}\geq2^{r+\alpha/2}d_{jk}$ for all $k$ (see (\ref{Relation})). Then,
\begin{eqnarray*}
&&
\frac{1}{\sqrt{2\uppi}}\int_{|t|\geq\delta_{j}/\sigma_{j}}\biggl( t-d_{jk}/\sigma_{j}\biggr) ^{2}
\mathrm{e}^{-(1/2)( t-d_{jk}/\sigma_{j}) ^{2}}\,\mathrm{d}t
\\
&&\quad\geq
\min_{0\leq x\leq2^{-(r+\alpha/2)}}\frac{1}{\sqrt{2\uppi}}\int_{|t|\geq 1}( t-x)^{2}
\mathrm{e}^{-(1/2)( t-x) ^{2}}\,\mathrm{d}t
\\
&&\quad\geq
\min_{0\leq x\leq2^{-1}}\frac{1}{\sqrt{2\uppi}}\int_{|t|\geq1}( t-x) ^{2}\mathrm{e}^{-(1/2)( t-x) ^{2}}\,\mathrm{d}t>0.57.
\end{eqnarray*}
Also, note that
\[
\frac{1}{\sqrt{2\uppi}}\int_{|t|<\delta_{j}/\sigma_{j}}\mathrm{e}^{-(1/2)( t-d_{jk}/\sigma_{j}) ^{2}}\,\mathrm{d}t\geq0.
\]
These two inequalities and (\ref{hatlambda}) imply that for all
$j>\hat{q}$,
\begin{eqnarray*}
\inf_{\delta_{j}\leq\sigma_{j}}\hat{\lambda}_{jk}&=&\inf_{\delta_{j}\leq
\sigma_{j}}\biggl\{ \sigma_{j}^{2}\frac{1}{\sqrt{2\uppi}}\int_{|t|>\delta
_{j}/\sigma_{j}}\biggl( t-\frac{d_{jk}}{\sigma_{j}}\biggr) ^{2}\mathrm{e}^{-(1/2)( t-d_{jk}/\sigma_{j}) ^{2}}\,\mathrm{d}t
\\
&&\hphantom{\inf_{\delta_{j}\leq\sigma_{j}}\biggl\{}{}+
d_{jk}^{2}\frac{1}{\sqrt{2\uppi}}\int_{|t|<\delta_{j}/\sigma_{j}}\mathrm{e}^{-(1/2)( t-d_{jk}/\sigma_{j}) ^{2}}\,\mathrm{d}t\biggr\}
\\
&\geq&
\inf_{\delta_{j}\leq\sigma_{j}}\biggl\{ \sigma_{j}^{2}\frac{1}{\sqrt{2\uppi
}}\int_{|t|>\delta_{j}/\sigma_{j}}\biggl( t-\frac{d_{jk}}{\sigma_{j}}\biggr)
^{2}\mathrm{e}^{-(1/2)( t-d_{jk}/\sigma_{j}) ^{2}}\,\mathrm{d}t\biggr\}
\geq0.57\sigma_{j}^{2}.
\end{eqnarray*}

For the case where {$\delta_{j}>\sigma_{j}$, we need some auxiliary
results.
Without loss of generality, we  let $d_{jk}\geq0$. First, note that if $%
\delta_{j}/\sigma_{j}>(1+\frac{d_{jk}}{\sigma_{j}})$, then}
\begin{eqnarray*}
&&
\frac{1}{\sqrt{2\uppi}}\int_{|t|<\delta_{j}/\sigma_{j}}\biggl[ \biggl( t-\frac{d_{jk}}{\sigma_{j}}\biggr) ^{2}-1\biggr]
\textrm{e}^{-(1/2)( t-d_{jk}/\sigma_{j}) ^{2}}\,\mathrm{d}t
\\
&&\quad\leq
\frac{1}{\sqrt{2\uppi}}\int_{-\infty}^{\infty}\biggl[ \biggl( t-\frac {d_{jk}}{\sigma_{j}}\biggr) ^{2}-1\biggr]
\textrm{e}^{-(1/2)( t-d_{jk}/\sigma_{j}) ^{2}}\,\mathrm{d}t=0
\end{eqnarray*}
so that
\[
\frac{1}{\sqrt{2\uppi}}\int_{|t|<\delta_{j}/\sigma_{j}}\biggl[ \biggl( t-\frac{%
d_{jk}}{\sigma_{j}}\biggr) ^{2}-1\biggr] \mathrm{e}^{-(1/2)( t-d_{jk}/\sigma_{j}) ^{2}}\,\mathrm{d}t\leq0
\]
and
\[
\frac{1}{\sqrt{2\uppi}}\int_{|t|<\delta_{j}/\sigma_{j}}\biggl( t-\frac{d_{jk}}{%
\sigma_{j}}\biggr) ^{2}\mathrm{e}^{-(1/2)( t-d_{jk}/\sigma_{j}) ^{2}}\,\mathrm{d}t\leq\frac{1}{\sqrt{2\uppi}}\int_{|t|<\delta_{j}/%
\sigma_{j}}\mathrm{e}^{-(1/2)( t-d_{jk}/\sigma_{j})^{2}}\,\mathrm{d}t.
\]
Similarly, if
$1\leq\delta_{j}/\sigma_{j}\leq(1+d_{jk}/\sigma_{j})$, then
\begin{eqnarray*}
&&
\frac{1}{\sqrt{2\uppi}}\int_{|t|<\delta_{j}/\sigma_{j}}\biggl[ \biggl( t-\frac{d_{jk}}{\sigma_{j}}\biggr) ^{2}-1\biggr]
\mathrm{e}^{-(1/2)( t-d_{jk}/\sigma_{j}) ^{2}}\,\mathrm{d}t
\\
&&\quad \leq
\frac{1}{\sqrt{2\uppi}}\int_{-\infty}^{\delta_{j}/\sigma_{j}}\biggl[\biggl( t-\frac{d_{jk}}{\sigma_{j}}\biggr) ^{2}-1\biggr]
\mathrm{e}^{-(1/2)( t-d_{jk}/\sigma_{j}) ^{2}}\,\mathrm{d}t.
\end{eqnarray*}
Moreover, since (\ref{Relation}), we have $%
d_{jk}/\sigma_{j}<1\leq\delta_{j}/\sigma_{j}$ so that an upper bound
is
given by%
\begin{eqnarray*}
&&
\frac{1}{\sqrt{2\uppi}}\int_{-\infty}^{d_{jk}/\sigma_{j}}\biggl[ \biggl( t-\frac{d_{jk}}{\sigma_{j}}\biggr) ^{2}-1\biggr]
\mathrm{e}^{-(1/2)( t-d_{jk}/\sigma_{j}) ^{2}}\,\mathrm{d}t
\\
&&\quad=
\frac{1}{\sqrt{2\uppi}}\int_{-\infty}^{0}[ t^{2}-1]
\mathrm{e}^{-(1/2)t^{2}}\,\mathrm{d}t=0.
\end{eqnarray*}
Hence, if $\delta_{j}>\sigma_{j}$, we also have the inequality
\[
\frac{1}{\sqrt{2\uppi}}\int_{|t|<\delta_{j}/\sigma_{j}}\biggl( t-\frac{d_{jk}}{%
\sigma_{j}}\biggr) ^{2}\mathrm{e}^{-(1/2)( t-d_{jk}/\sigma_{j}
) ^{2}}\,\mathrm{d}t\leq\frac{1}{\sqrt{2\uppi}}\int_{|t|<\delta_{j}/%
\sigma_{j}}\mathrm{e}^{-(1/2)( t-d_{jk}/\sigma_{j})
^{2}}\,\mathrm{d}t.
\]
In summary, we obtain
\begin{eqnarray*}
&&
\frac{1}{\sqrt{2\uppi}}\int_{|t|>\delta_{j}/\sigma_{j}}\biggl( t-\frac{d_{jk}}{%
\sigma_{j}}\biggr) ^{2}\mathrm{e}^{-(1/2)( t-d_{jk}/\sigma_{j}
) ^{2}}\,\mathrm{d}t
\\
&&\quad=
1-\frac{1}{\sqrt{2\uppi}}\int_{|t|<\delta_{j}/\sigma_{j}}\biggl( t-\frac {%
d_{jk}}{\sigma_{j}}\biggr) ^{2}\mathrm{e}^{-(1/2)( t-d_{jk}/\sigma_{j}) ^{2}}\,\mathrm{d}t
\\
&&\quad\geq
1-\frac{1}{\sqrt{2\uppi}}\int_{|t|<\delta_{j}/\sigma_{j}}\mathrm{e}^{-(1/2)( t-d_{jk}/\sigma_{j}) ^{2}}\,\mathrm{d}t.
\end{eqnarray*}
Defining%
\[
\gamma=\frac{1}{\sqrt{2\uppi}}\int_{|t|<\delta_{j}/\sigma_{j}}\mathrm{e}^{-(1/2)( t-d_{jk}/\sigma_{j}) ^{2}}\,\mathrm{d}t\in\lbrack0,1],
\]
this inequality, together with (\ref{hatlambda}), implies that for all
$j\geq
\hat{q}$, and all $k$ and $n$ large enough, $\inf_{\delta_{j}>\sigma_{j}}\hat{%
\lambda}_{jk}$ is equal to%
\begin{eqnarray*}
&&\inf_{\delta_{j}>\sigma_{j}}\biggl\{ \sigma_{j}^{2}\frac{1}{\sqrt{2\uppi}}%
\int_{|t|>\delta_{j}/\sigma_{j}}\biggl( t-\frac{d_{jk}}{\sigma_{j}}\biggr)
^{2}\mathrm{e}^{-(1/2)( t-d_{jk}/\sigma_{j})
^{2}}\,\mathrm{d}t
\\
&&\hphantom{\inf_{\delta_{j}>\sigma_{j}}\biggl\{}{}+d_{jk}^{2}\frac{1}{\sqrt{2\uppi}}\int_{|t|<\delta_{j}/\sigma_{j}}\mathrm{e}^{-(1/2)( t-d_{jk}/\sigma_{j}) ^{2}}\,\mathrm{d}t\biggr\}
\end{eqnarray*}
so that
\[
\inf_{\delta_{j}>\sigma_{j}}\hat{\lambda}_{jk}\geq\inf_{\gamma\in\lbrack
0,1]}\{ (1-\gamma)\sigma_{j}^{2}+\gamma d_{jk}^{2}\} =d_{jk}^{2}.
\]
Moreover, note that the minimum is attained at the border. Now,
\begin{eqnarray*}
\inf_{\delta_{j}}\hat{\lambda}_{jk}&=&\min\Bigl\{ \inf_{\delta_{j}\leq
\sigma_{j}}\hat{\lambda}_{jk},\inf_{\delta_{j}>\sigma_{j}}\hat{\lambda}%
_{jk}\Bigr\} \geq\min\{ 0.57\sigma_{j}^{2},d_{jk}^{2}\}
\\
&\geq&\min\{ 0.57\cdot2^{2r+\alpha}\cdot\max
d_{jk}^{2},d_{jk}^{2}\} =d_{jk}^{2}, %,
\end{eqnarray*}
where the last inequality follows from (\ref{Relation}). Clearly,
the value of $d_{jk}^{2}$ is attained if and only if
$\delta_{j}=\infty$.

Finally, we obtain
\begin{eqnarray*}
\min_{\delta_{j}}\lambda_{jk}
&=&\min_{\delta_{j}}\hat{\lambda}_{jk}+\mathrm{O}\bigl(
2^{\alpha(J+j)/2}n^{-(1+\alpha/2)}\bigr) +\mathrm{O}\bigl( n^{-(1-\alpha
/2)}2^{-(2r+\alpha/2)(J+j)}\bigr)
\\
&=&d_{jk}^{2}+\mathrm{O}\bigl( 2^{\alpha(J+j)/2}n^{-(1+\alpha/2)}\bigr)+\mathrm{O}\bigl( n^{-(1-\alpha
/2)}2^{-(2r+\alpha/2)(J+j)}\bigr).
\end{eqnarray*}

Now, $d_{jk}=\mathrm{O}(2^{-((2r+1)/2)(J+j)})$ and the assumption
\[
{\hat{q}< j <
\frac{2+\alpha}{4r+2+\alpha}\log_2n-J}
\]
implies that
\[
2^{-(2r+1)(J+j)}>2^{\alpha(J+j)/2}n^{-(1+\alpha/2)}
\]
and
\[
2^{\alpha(J+j)/2}n^{-(1+\alpha/2)}> n^{-(1-\alpha /2)}2^{-(2r+\alpha/2)(J+j)}.
\]
Therefore, the remainder term $2^{\alpha(J+j)/2}n^{-(1+\alpha/2)}$ is
of smaller order than $d_{jk}^2$, and
$\mathrm{O}(2^{\alpha(J+j)/2}\times n^{-(1+\alpha/2)})$ dominates $\mathrm{O}(n^{-(1-\alpha
/2)}2^{-(2r+\alpha/2)(J+j)})$. This completes the proof of Lemma
\ref{Lambda_6}.

We now come back to the proof of Theorem \ref{theorema1}. Suppose that $\phi$ and $\psi$ are piecewise
differentiable. We define
\[
\hat{J}=\log_{2}n^{\alpha/(2r+\alpha)}+\frac{1}{2r+\alpha}\log
_{2}\biggl( \frac{\nu_{r}^{2}}{(r!)^{2}C_{\psi}^{2}N^{2r+\alpha}}\max _{t\in
[ 0,1]}\bigl[ g^{(r)}(t)\bigr] ^{2}\biggr)+1
\]
and  let $J \ge \hat{J}$. Noting that $\Lambda_{i} \ge 0$
($i=1,2,3,4$) and taking into account Lemma \ref{Lambda_2}, we obtain,
for all $q \ge 0,$
\begin{eqnarray*}
E\int_{0}^{1}\bigl( g(t)-\hat {g}%
(t)\bigr) ^{2}\,\mathrm{d}t&=& \Lambda_{1}+\Lambda_{2}+\Lambda_{3}+\Lambda_{4}
 \ge \Lambda_{2}
 \\
&\ge& C_{\phi}^{2}n^{-\alpha}N^{\alpha }2^{\alpha J}+
\mathrm{O}(n^{-1}2^{J})+\mathrm{O}\bigl(n^{-\alpha}2^{-J(1-\alpha)}\bigr) \ge
C_1n^{-2r\alpha/(2r+\alpha)}.
\end{eqnarray*}

Now, consider $J< \hat J$ and let $q \le \hat{q}$, where $\hat{q}$ is
defined as in Lemma \ref{Lambda_6}. Lemmas \ref{Lambda_4} and~\ref{Lambda_5} imply that
\[
\sum_{j=\hat{q}+1}^{\infty}\sum _{k=-N+1}^{N2^{J+j}-1}d_{jk}^{2}\ge
C_2n^{-2r\alpha/(2r+\alpha)}.
\]
Since $q \le \hat{q}$, we have
\begin{eqnarray*}
E\int_{0}^{1}\bigl( g(t)-\hat {g}%
(t)\bigr) ^{2}\,\mathrm{d}t&=& \Lambda_{1}+\Lambda_{2}+\Lambda_{3}+\Lambda_{4}
\geq\Lambda _{4} = \sum_{j=q+1}^{\infty}\sum
_{k=-N+1}^{N2^{J+j}-1}d_{jk}^{2}\geq\sum_{j=\hat{q}+1}^{\infty}\sum
_{k=-N+1}^{N2^{J+j}-1}d_{jk}^{2}\\
&\geq&
C_{2}n^{-2r\alpha/(2r+\alpha)}.
\end{eqnarray*}

For the other case, where $q>\hat{q}$, taking into account
$\Lambda_{3}$ in (\ref{MISE}) and Lemma \ref{Lambda_6} leads to
\begin{eqnarray*}
&&E\int_{0}^{1}\bigl( g(t)-\hat {g}%
(t)\bigr) ^{2}\,\mathrm{d}t
\\
&&\quad\ge \Lambda_{3}=\sum_{j=0}^{q}\sum_{k=-N+1}^{N2^{J+j}-1}\{ E[ (\hat {d}%
_{jk}-d_{jk})^{2}I(|\hat{d}_{jk}|>\delta_{j})] +E[ d_{jk}^{2}I(|%
\hat{d}_{jk}|\leq\delta_{j})] \}
\\
&&\quad=\sum_{j=0}^{q}\sum_{k=-N+1}^{N2^{J+j}-1}\lambda_{jk} \ge \sum_{j=0}^{\hat{q}}\sum_{k=-N+1}^{N2^{J+j}-1}\lambda_{jk}+\sum_{j=\hat{q}%
+1}^{q}\sum_{k=-N+1}^{N2^{J+j}-1}\min_{\delta_{j}}\lambda_{jk}
\\
&&\quad\ge \sum_{j=\hat{q}%
+1}^{\hat{q}+1}\sum_{k=-N+1}^{N2^{J+j}-1}\min_{\delta_{j}}\lambda_{jk}=\sum_{j=\hat{q}%
+1}^{\hat{q}+1}\sum_{k=-N+1}^{N2^{J+j}-1}d_{jk}^2
+\mathrm{O}\bigl( 2^{(1+\alpha/2)(J+\hat q)}n^{-(1+\alpha/2)}\bigr)
\\
&&\quad\ge
C_3n^{-2r\alpha/(2r+\alpha)}.
\end{eqnarray*}
In summary, we have obtained a lower bound:
\[
\min_{\{ \delta_{j}\} ,q,J}E\int_{0}^{1}\bigl( g(t)-\hat {g}%
(t)\bigr) ^{2}\,\mathrm{d}t=\min_{\{ \delta_{j}\} ,q,J}(
\Lambda_{1}+\Lambda_{2}+\Lambda_{3}+\Lambda_{4}) \ge
Cn^{-2r\alpha/(2r+\alpha)}.
\]
It is shown in the proof of Theorem \ref{Theorem2} that equality can
indeed be achieved, by using a~specific choice of $\delta_{j}$, $q$,
$J$ and $C$. This completes the proof of Theorem \ref{theorema1}.
\end{pf}

\begin{pf*}{Proof of Theorem \ref{Theorem2}}
Under the conditions of Theorem \ref{Theorem2}, and taking into account
Lemmas \ref{Lambda_3} and \ref{Lambda_4}, we obtain
that
\begin{eqnarray*}
\Lambda_{3}&=&\sum_{j=0}^{q}\sum_{k=-N+1}^{N2^{J+j}-1}E\{ [ \hat {d}%
_{jk}-d_{jk}] ^{2}\} =\sum_{j=0}^{q}\sum
_{k=-N+1}^{N2^{J+j}-1}\bigl(\sigma_{j}^{2}+\bigl( E( \hat d_{jk})-d_{jk}\bigr) ^{2}\bigr)
\\[-2pt]
&=&\frac{C_{\psi}^{2}}{2^{\alpha}-1}\bigl( 2^{\alpha(q+1)}-1\bigr) N^{\alpha
}n^{-\alpha}2^{\alpha
J}+\mathrm{O}(n^{-1}2^{J+q})+\mathrm{O}\bigl(n^{-\alpha}2^{-J(1-\alpha)}\bigr)+\mathrm{O}\bigl(n^{-2}2^{2(J+q)}\bigr).
\end{eqnarray*}
This, together with Lemmas \ref{Lambda_1}, \ref{Lambda_2} and \ref{Lambda_5},
implies that the expression in (\ref{MISE}) will take the following form:
\begin{eqnarray}\label{MISEend}
\mathit{MISE}_{g}(q,J)&=&\Lambda_{1}+\Lambda_{2}+\Lambda_{3}+\Lambda_{4}\nonumber
\\[-2pt]
&=&\biggl( C_{\phi}^{2}-\frac{C_{\psi}^{2}}{2^{\alpha}-1}\biggr) N^{\alpha
}n^{-\alpha}2^{\alpha J}+\frac{2^{\alpha}C_{\psi}^{2}}{2^{\alpha}-1}%
N^{\alpha }n^{-\alpha}2^{\alpha(J+q)}\nonumber
\\[-9pt]\\[-9pt]
&&{}+\frac{1}{(r!)^{2}}\frac{1}{2^{2r}-1}N^{-2r}2^{-2r(J+q)}\int_{0}^{1}%
\nu_{r}^{2}\bigl(g^{(r)}(t)\bigr)^2\,\mathrm{d}t\nonumber
\\[-2pt]
&&{}+\mathrm{o}\bigl( 2^{-2r(J+q)}\bigr) +\mathrm{O}(n^{-1}2^{J+q})+\mathrm{O}\bigl(n^{-\alpha}2^{-J(1-\alpha)}\bigr).\nonumber
\end{eqnarray}
Now, let $q$ and $J$ be such that{ $\mathit{MISE}$} is minimal. Then,
by (\ref{MISE}), $\delta_j=0$ and
\[
\mathit{MISE}_{g}(q,J)-\mathit{MISE}_{g}(q+1,J)<0
\]
imply that
\[
\mathit{MISE}_{g}(q,J)-\mathit{MISE}_{g}(q+1,J)=
\sum_{k=-N+1}^{N2^{J+q+1}-1}d_{q+1,k}^{2}-\sum_{k=-N+1}^{N2^{J+q+1}-1}%
\sigma_{q+1}^{2}<0.
\]
By an argument analogous to the one used in the proof of (\ref{Relation}), the last
inequality, together with Lemmas \ref{Lambda_3} and
\ref{Lambda_4}, implies that for $n$ large enough, we have
\[
C_{\psi}^{2}n^{-\alpha}N^{\alpha}2^{\alpha(J+q+1)}\ge\frac{1}{(r!)^{2}}%
N^{-2r}2^{-2r(J+q+1)}\int_{0}^{1}\nu_{r}^{2}\bigl(g^{(r)}(t)\bigr)^2\,\mathrm{d}t
\]
and
\begin{equation}\label{Inequaliy1}
q\ge\log_{2}n^{\alpha/(2r+\alpha)}-J-1+\frac{1}{2r+\alpha}\log_{2}\biggl[
\frac{\int_{0}^{1}\nu_{r}^{2}(g^{(r)}(t))^2\,\mathrm{d}t}{C_{\psi}^{2}(r!)^{2}}\biggr]
-\log_{2}N.
\end{equation}
On the other hand,
\[
\mathit{MISE}_{g}(q,J)-\mathit{MISE}_{g}(q-1,J)<0
\]
implies the second necessary condition
\[
C_{\psi}^{2}n^{-\alpha}N^{\alpha}2^{\alpha(J+q)}\le\frac{1}{(r!)^{2}}%
N^{-2r}2^{-2r(J+q)}\int_{0}^{1}\nu_{r}^{2}\bigl(g^{(r)}(t)\bigr)^2\,\mathrm{d}t\vadjust{\goodbreak}
\]
so that
\begin{equation}\label{Inequaliy2}
q\le\log_{2}n^{\alpha/(2r+\alpha)}+\frac{1}{2r+\alpha}\log_{2}\biggl[
\frac{\int_{0}^{1}\nu_{r}^{2}(g^{(r)}(t))^2\,\mathrm{d}t}{C_{\psi}^{2}(r!)^{2}}\biggr]
-\log_{2}N-J.
\end{equation}
Note that $q$ and $J$ are integers. The inequalities
(\ref{Inequaliy1}) and (\ref{Inequaliy2}) then imply that the value
\begin{equation}\label{bestq}
q^*=\biggl\lfloor\log_{2}n^{\alpha/(2r+\alpha)}+\frac{1}{2r+\alpha}\log_{2}\biggl[
\frac{\int_{0}^{1}\nu_{r}^{2}(g^{(r)}(t))^2\,\mathrm{d}t}{C_{\psi}^{2}(r!)^{2}}\biggr]
-\log_{2}N\biggr\rfloor-J
\end{equation}
asymptotically minimizes the $\mathit{MISE}$. Using the definition of
$\Delta_n(g,C_{\psi})$ in (\ref{DeltaPsi}), we conclude that
\[
q^*=\log_{2}n^{\alpha/(2r+\alpha)}+\frac{1}{2r+\alpha}\log_{2}\biggl[
\frac{\int_{0}^{1}\nu_{r}^{2}(g^{(r)}(t))^2\,\mathrm{d}t}{C_{\psi}^{2}(r!)^{2}}\biggr]
-\log_{2}N-J-\Delta_n(g,C_{\psi}).
\]
Note that if $\Delta_n(g,C_{\psi})\ne 0$, then for every fixed $J,$
there exists a unique $q^*$ such that (\ref{Inequaliy1}) and
(\ref{Inequaliy2}) hold.

Combining these results with (\ref{MISEend}) yields
\begin{eqnarray}\label{bestMISEJ}
\mathit{MISE}_{g}(q^*,J)
&=&
2^{-\alpha \Delta_n(g,C_{\psi})}\biggl(C_{\phi}^{2}-\frac{C_{\psi}^{2}}{2^{\alpha}-1}\biggr)N^{\alpha}n^{-\alpha}2^{\alpha J}\nonumber
\\[-2pt]
&&{}+
\biggl( \frac{2^{2r\Delta_n(g,C_{\psi})}}{2^{2r}-1}+\frac{2^{\alpha(1-\Delta_n(g,C_{\psi}))}}{2^{\alpha}-1}\biggr)\nonumber
\\[-9pt]\\[-9pt]
&&\hphantom{{}+{}}{}\times
C_{\psi }^{4r/(2r+\alpha)}\biggl( \frac{\nu_{r}^{2}}{(r!)^{2}}\int_{0}^{1}\bigl(g^{(r)}(t)\bigr)^2\,\mathrm{d}t\biggr) ^{\alpha/(2r+\alpha)}n^{-2r\alpha /(2r+\alpha)}\nonumber
\\[-2pt]
&&{}+\mathrm{O}(n^{-1}2^{J})+\mathrm{o}\bigl(n^{-2r\alpha/(2r+\alpha)}\bigr)+\mathrm{O}\bigl(n^{-\alpha}2^{-J(1-\alpha)}\bigr).\nonumber
\end{eqnarray}
The first term is monotonically decreasing in $J$ if
$(2^{\alpha}-1)C_{\phi}^{2}<C_{\psi}^{2}$, and monotonically increasing
if $(2^{\alpha}-1)C_{\phi}^{2}>C_{\psi}^{2}$. The second term does not
depend of $J$. Hence, if $(2^{\alpha}-1)C_{\phi}^{2}>C_{\psi}^{2}$,
then the optimal decomposition level $J$ is equal to zero. Note that
the optimal decomposition level is not unique since the same asymptotic
expressions will be achieved for all integers $J$ such that
$2^J=\mathrm{o}(n^{\alpha/(2r+\alpha)})$. Combining this with the previous
formulas implies that $\mathit{MISE}_{g}(q,J)$ is equal to
\begin{eqnarray}\label{BestMISE}
&&\hspace*{-15pt}\biggl( \frac{2^{2r\Delta_n(g,C_{\psi})}}{2^{2r}-1}+
\frac{2^{\alpha(1-\Delta_n(g,C_{\psi}))}}{2^{\alpha}-1}\biggr)C_{\psi }^{4r/(2r+\alpha)}\biggl( \frac{\nu_{r}^{2}}{(r!)^{2}}%
\int_{0}^{1}\bigl(g^{(r)}(t)\bigr)^2\,\mathrm{d}t\biggr) ^{\alpha/(2r+\alpha)}n^{-2r\alpha/(2r+\alpha)}\nonumber\hspace*{15pt}
\\[-9pt]\\[-9pt]
&&\hspace*{-15pt}\quad{}+
\mathrm{o}\bigl( n^{-2r\alpha/(2r+\alpha)}\bigr) .\hspace*{15pt}\nonumber
\end{eqnarray}

On the other hand, suppose that
$(2^{\alpha}-1)C_{\phi}^{2}<C_{\psi}^{2}$. Taking into account
(\ref{bestq}) and $q\ge0$ (see (\ref{4})), we then have
\[
0 \le J \le
\biggl\lfloor\log_{2}n^{\alpha/(2r+\alpha)}+\frac{1}{2r+\alpha}\log_{2}\biggl[
\frac{\int_{0}^{1}\nu_{r}^{2}(g^{(r)}(t))^2\,\mathrm{d}t}{C_{\psi}^{2}
(r!)^{2}}\biggr] -\log_{2}N\biggr\rfloor.\vadjust{\goodbreak}
\]
Hence, the optimal choice of $J$ is
\begin{equation}\label{Temp}
J=\biggl\lfloor\log_{2}n^{\alpha/(2r+\alpha)}+\frac{1}{2r+\alpha}\log_{2}\biggl[
\frac{\int_{0}^{1}\nu_{r}^{2}(g^{(r)}(t))^2\,\mathrm{d}t}{C_{\psi}^{2}
(r!)^{2}}\biggr] -\log_{2}N\biggr\rfloor.
\end{equation}
Due to (\ref{bestq}), this also implies that $q^*=0$.

Note that (\ref{4}) with $q \ge 0$ and $\delta_j=0$ always includes
at least one level of mother wavelets. The case where the estimate
includes father wavelets only is automatically considered in Theorem
\ref{theorema1}, namely, if $q=0$ and $\delta_0=\infty$. To complete
the proof, we also need to compare with the estimate that only includes
father wavelets. Thus, we consider
\[
\tilde{g}(t)=\sum_{k=-N+1}^{N2^{J}-1}\hat{s}_{Jk}\phi_{Jk}(t)
\]
and denote the corresponding mean integrated square error by
$\mathit{MISE}_g(-1,J)$. Then,
\begin{eqnarray*}
\mathit{MISE}_g(-1,J)&=&\sum_{k=-N+1}^{N2^{J}-1} [ E( \hat{s}_{Jk})
-s_{Jk}] ^{2}+\sum _{k=-N+1}^{N2^{J}-1}E\{ [
\hat{s}_{Jk}-E(\hat{s}_{Jk})] ^{2}\}
\\
&&{}+\sum_{j=0}^{\infty}\sum_{k=-N+1}^{N2^{J+j}-1}d_{jk}^{2}.
\end{eqnarray*}
Let $J^*$ be such that{ $\mathit{MISE}_g(-1,J^*)$} is minimal. Then,
\[
\mathit{MISE}_g(-1,J^*)-\mathit{MISE}_g(-1,J^*+1) < 0
\]
and
\[
\mathit{MISE}_g(-1,J^*)-\mathit{MISE}_g(-1,J^*-1) < 0.
\]
Suppose that $n$ is large enough. Elementary calculations
similar to those above then show that the optimal decomposition level
$J^*$ is given by
\begin{equation}\label{bestJ}
J^*=\biggl\lfloor\log_{2}n^{\alpha/(2r+\alpha)}+\frac{1}{2r+\alpha}\log_{2}\biggl[
\frac{\int_{0}^{1}\nu_{r}^{2}(g^{(r)}(t))^2\,\mathrm{d}t}{C_{\phi}^{2}(2^{\alpha
}-1)(r!)^{2}}\biggr] -\log_{2}N\biggr\rfloor+1.
\end{equation}
Defining $\Delta_n(g,C_{\phi})$ as in (\ref{DeltaPhi}), the
corresponding $\mathit{MISE}$ is equal to
\begin{eqnarray*}
&&\biggl( \frac{2^{2r\Delta_n(g,C_{\phi})}}{2^{2r}-1}+
\frac{2^{\alpha(1-\Delta_n(g,C_{\phi}))}}{2^{\alpha}-1}\biggr) \bigl(
C_{\phi}^{2}(2^{\alpha}-1)\bigr) ^{2r/(2r+\alpha)}
\\
&&\quad{}\times\biggl( \frac {%
\nu_{r}^{2}}{(r!)^{2}}\int_{0}^{1}\bigl(g^{(r)}(t)\bigr)^2\,\mathrm{d}t\biggr) ^{\alpha/(2r+\alpha)}n^{-2r\alpha/(2r+\alpha)}
+\mathrm{o}\bigl( n^{-2r\alpha/(2r+\alpha)}\bigr).
\end{eqnarray*}
Now, let $(2^{\alpha}-1)C_{\phi}^{2}>C_{\psi}^{2}$. Suppose that $J$
defined by (\ref{bestJ}) and the estimator consisting of only father
wavelets minimizes the MISE. Now,
\begin{eqnarray*}
&&\mathit{MISE}_g(0,J)-\mathit{MISE}_g(-1,J+1)
\\
&&\quad=n^{-\alpha}N^{\alpha}2^{\alpha
J}\bigl(C_{\psi}^2-C_{\phi}^2(2^{\alpha}-1)\bigr)+\mathrm{o}\bigl(n^{-2r\alpha/(2r+\alpha)}\bigr)
\end{eqnarray*}
so that, for $n$ large enough,
\[
\mathit{MISE}_g(0,J)-\mathit{MISE}_g(-1,J+1)<0,
\]
which is a contradiction. It thus follows that the best $J$ is
equal to zero, $q$ is defined by~(\ref{bestq}) and the MISE is as
in (\ref{BestMISE}).

Now, suppose that
\[
C_{\phi}^2(2^{\alpha}-1)<C_{\psi}^2,
\]
$q=0$ and $J$ given by (\ref{Temp}) minimizes the MISE. Consider
\begin{eqnarray*}
&&\mathit{MISE}_g(-1,J+1)-\mathit{MISE}_g(0,J)
\\
&&\quad=n^{-\alpha}N^{\alpha}2^{\alpha
J}\bigl(C_{\phi}^2(2^{\alpha}-1)-C_{\psi}^2\bigr)+\mathrm{o}\bigl(n^{-2r\alpha/(2r+\alpha)}\bigr).
\end{eqnarray*}
Using the same argument as before,
$\mathit{MISE}_g(-1,J+1)-\mathit{MISE}_g(0,J)<0$ for $n$ large enough.
Thus, the best estimator includes only father wavelets and the optimal
decomposition level is defined by (\ref{bestJ}).

In conclusion, we consider the case $\Delta_n(g,C_{\psi})=0$. Suppose
that
\[
(2^{\alpha}-1)C_{\phi}^{2}>C_{\psi}^{2},
\]
$J=0$ and $q$ as in (\ref{bestq}) minimizes the MISE. Now,
\begin{eqnarray*}
&&\mathit{MISE}_g(q,0)-\mathit{MISE}_g(q-1,0)
\\
&&\quad=\biggl( \frac{1}{2^{2r}-1}+ \frac{2^{\alpha}}{2^{\alpha}-1}-
\frac{2^{2r}}{2^{2r}-1}
-\frac{1}{2^{\alpha}-1}\biggr)
\\
&&\qquad{}\times C_{\psi }^{4r/(2r+\alpha)}\biggl( \frac{\nu_{r}^{2}}{(r!)^{2}}%
\int_{0}^{1}\bigl(g^{(r)}(t)\bigr)^2\,\mathrm{d}t\biggr) ^{\alpha/(2r+\alpha)}n^{-2r\alpha/(2r+\alpha)}
\\
&&\qquad{}+\mathrm{o}\bigl( n^{-2r\alpha/(2r+\alpha)}\bigr)=\mathrm{o}\bigl(
n^{-2r\alpha/(2r+\alpha)}\bigr).
\end{eqnarray*}
Then, for every fixed $J,$ there exist two smoothing parameters that
minimize the MISE asymptotically. The same also follows  for the
case $(2^{\alpha}-1)C_{\phi}^{2}<C_{\psi}^{2}$ and
$\Delta_n(g,C_{\phi})=0$.  This completes the proof.
\end{pf*}

\begin{pf*}{Proof of Theorem \ref{theorema3}}
The extension to functions with
piecewise continuous $r$th derivatives follows from the following
lemma, which can be proven in a similar manner as
Lemmas 4.5 and 4.6 in Li and Xiao \cite{r39}.
\end{pf*}

\begin{lemma}
\label{Lambda_7}
Suppose that the assumptions of Theorem
\ref{theorema3} hold. Then:
\begin{longlist}
\item[(i)] if $(2^{\alpha}-1)C_{\phi}^{2}>C_{\psi}^{2}$, then
\begin{eqnarray*}
&&\sum_{j=q^*+1}^{q}\sum_{k=-N+1}^{N2^{J+j}-1}\lambda_{jk}
+\sum_{j=q+1}^{\infty}\sum_{k=-N+1}^{N2^{J+j}-1}d_{jk}^{2}
\\
&&\quad=\frac{\nu_{r}^{2}}{(r!)^{2}}\frac{1}{2^{2r}-1}%
N^{-2r}2^{-2r(J+q^*)}\int_{0}^{1}\bigl(g^{(r)}(t)\bigr)^2\,\mathrm{d}t+\mathrm{o}\bigl(
2^{-2r(J+q)}\bigr) ;
\end{eqnarray*}
\item[(ii)] if $(2^{\alpha}-1)C_{\phi}^{2}<C_{\psi}^{2}$, then
\begin{eqnarray*}
&&\sum_{j=0}^{q}\sum_{k=-N+1}^{N2^{J+j}-1}\lambda_{jk}
+\sum_{j=q+1}^{\infty}\sum_{k=-N+1}^{N2^{J+j}-1}d_{jk}^{2}
\\
&&\quad=\frac{\nu_{r}^{2}}{(r!)^{2}}\frac{1}{2^{2r}-1}%
N^{-2r}2^{-2r(J-1)}\int_{0}^{1}\bigl(g^{(r)}(t)\bigr)^2\,\mathrm{d}t+\mathrm{o}(
2^{-2rJ}).
\end{eqnarray*}
\end{longlist}
\end{lemma}
\end{appendix}

\section*{Acknowledgements}

This work was supported in part by a grant from the German
Research Foundation (Research Unit 518). We would like to thank the
referees for careful reading of previous versions of the paper and
for constructive suggestions that helped to improve the quality of the
presentation.

\printhistory

\end{document}